\documentclass[]{amsart}
\usepackage{color}

 \usepackage[latin1]{inputenc}
 \usepackage[T1]{fontenc}
 \usepackage[normalem]{ulem}
\usepackage{verbatim}
 \usepackage{graphicx}
 \usepackage{epstopdf}

\usepackage[latin1]{inputenc}
\usepackage{amsmath}
\usepackage{amssymb}
\usepackage{amsfonts}
\usepackage{amsthm}
\usepackage{bbm}

\newcommand{\C}{\mathbb{C}}
\newcommand{\N}{\mathbb{N}}
\newcommand{\R}{\mathbb{R}}

\newtheorem{theorem}{Theorem}[section]

\theoremstyle{remark}

\theoremstyle{definition}

\newtheorem{conjecture}[theorem]{Conjecture}
\newtheorem*{merci}{Acknowledgements}

\begin{document}

\title[Davey-Stewartson]{Numerical study
of blow-up mechanisms for Davey-Stewartson II systems}
\author{Christian Klein}
\address{Institut de Math\'ematiques de Bourgogne, UMR 5584\\
                Universit\'e de Bourgogne-Franche-Comt\'e, 9 avenue Alain Savary, 21078 Dijon
                Cedex, France\\
    E-mail Christian.Klein@u-bourgogne.fr}

\author{Nikola Stoilov}
\address{Institut de Math\'ematiques de Bourgogne, UMR 5584\\
                Universit\'e de Bourgogne-Franche-Comt\'e, 9 avenue Alain Savary, 21078 Dijon
                Cedex, France\\
    E-mail Nikola.Stoilov@u-bourgogne.fr}
\date{\today}
\maketitle
\large

\begin{abstract}
We present a detailed numerical study of various blow-up issues in 
the context of the focusing Davey-Stewartson II equation. To this end 
we study Gaussian initial data and perturbations of the lump and the 
explicit blow-up solution due to Ozawa. Based on the numerical 
results it is conjectured that the 
blow-up in all cases is self similar, and that the time dependent 
scaling is as in the Ozawa solution and not as in the stable blow-up 
of standard $L^{2}$ critical nonlinear Schr\"odinger equations. The 
blow-up profile is given by a 
dynamically rescaled lump.

\end{abstract}

\section{Introduction}
This paper is devoted to the study of blow-up mechanisms for 
solutions to the focusing
 Davey-Stewartson (DS) II  equation
\begin{equation}
    \label{DSII}
\begin{array}{ccc}
i
\partial_{t}\psi+\partial_{xx}\psi-\partial_{yy}\psi-2\left(\Phi+\left|\psi\right|^{2}\right)\psi & = & 0,\quad \psi :  \R^2\times \R \to \C,
\\
\partial_{xx}\Phi+\partial_{yy}\Phi+2\partial_{xx}\left|\psi\right|^{2} & = & 0, \quad \Phi :  \R^2\times \R \to \R,\\
\psi(.,0)=\psi_0.
\end{array}
\end{equation}
Equation (\ref{DSII}) is  a completely integrable \cite{AH,Sc} 
two-dimensional nonlinear Schr\"odinger (NLS) 
equation which appeared first in the context of water waves 
(\cite{DS,DR,AS}, see \cite{La} for more references and a rigorous 
justification). It also appears  in the 
context of ferromagnetism \cite{Le}, plasma physics \cite{MRZ} and 
nonlinear optics \cite{NM}.
Eliminating the \emph{mean field} $\Phi$ from (\ref{DSII}) 
by inverting the Laplace operator with some fall off condition at 
infinity, one gets as in \cite{KS}
\begin{equation}\label{box}
i \partial_{t}\psi+\Box\psi+2\lbrack(\Delta^{-1}\Box) |\psi|^2\rbrack\psi=0,
\end {equation}
where $\Box =\partial_{xx}-\partial_{yy}$, and where 
$\Delta=\partial_{xx}+\partial_{yy}$ are the D'Alembert and Laplace 
operator respectively. The DS II equation can thus 
be seen as a nonlocal hyperbolic NLS equation. 

The DS II equation as a special case of 2d NLS equations
has a well known scaling invariance: if 
$\psi(x,y,t)$ is a solution to (\ref{DSII}), so is $\lambda\psi(x/\lambda 
x,y/\lambda,t/\lambda^{2})$ with constant $\lambda\in 
\mathbb{R}/\{0\}$. Since the $L^{2}$ norm of $\psi$ is invariant 
under this rescaling, the equation is  $L^{2}$ critical.
It is well known that $L^{2}$ critical focusing NLS  equations 
can have solutions with a blow-up in finite time, see \cite{SS99} for 
a detailed discussion and for references. 
The above rescaling can be studied with a time dependent scaling 
factor $L(t)$ (which is not a symmetry of the equation),
\begin{equation}\label{resc}
    X = \frac{x}{L(t)},\quad Y = \frac{y}{L(t)},\quad 
    \tau=\int_{0}^{t}\frac{dt'}{L^{2}(t')},\quad \Psi(\xi,\eta,\tau) = 
    L(t)\psi(x,y,t).
\end{equation}
It was shown by Merle and Rapha\"el \cite{MeRa} that in cases with a 
generic blow-up in one point, the blow-up of the NLS solution
is self-similar and follows 
the `dynamical rescaling' (\ref{resc}) with a scaling factor of the 
form 
\begin{equation}
    L(t) \propto \sqrt{\frac{t~^{*}-t}{\ln|\ln(t^{*}-t)|}}
    \label{loglog},
\end{equation}
where $t^{*}$ is the blow-up time. 
If the spatial operators  in 
(\ref{DSII}) are both elliptic, the 
resulting equation is the elliptic-elliptic DS following the notation 
of
\cite{GS}. Numerically its solutions appear  to have finite time blow-up, see \cite{PSSW}, but a 
proof for the blow-up mechanism as in \cite{MeRa} 
for the cubic NLS has not yet been established. On 
the other hand a global well posedness result for the hyperbolic NLS equation 
$$i \partial_{t}\psi+\Box\psi+ |\psi|^2 \psi=0
$$
appears to be given in \cite{Totz}. Thus it 
is not obvious whether the focusing DS II equation should have solutions with a  
blow-up in finite time. But equation (\ref{DSII}) has an additional 
symmetry sometimes referred to as pseudoconformal invariance: if 
$\psi(x,y,t)$ is a solution to the DS II equation for $t>0$, so is 
\begin{equation}
    \tilde{\psi}(x,y,t) = \exp\left(\frac{i(x^{2}-y^{2})}{4t}\right)
\psi\left(\frac{x}{t},\frac{y}{t},\frac{1}{t}\right).
    \label{pseudo}
\end{equation}
This symmetry was used by Ozawa \cite{Oz} to construct an exact 
blow-up solution starting from the \emph{lump} solution (\ref{lump}) of the 
focusing DS II solution, a solitary wave solution with algebraic 
decay towards infinity. The blow-up in the Ozawa solution 
(\ref{soluoz}) is again self-similar, i.e., of the form 
(\ref{resc}), but this time with a 
scaling factor
\begin{equation}
    L(t)\propto t^{*}-t.
    \label{ozawascal}
\end{equation}
Note that a similar blow-up mechanism based on the pseudoconformal 
invariance is also known for the standard 
$L^{2}$ critical NLS equations, but that it is unstable in this case, see 
\cite{MeRa}. 

The question to be addressed numerically in this paper is whether 
the blow-up in DS II solutions observed in various examples is of the type (\ref{resc}), 
and if so whether it has the scaling (\ref{loglog}) of the generic 
blow-up in NLS or the scaling (\ref{ozawascal}) of the Ozawa 
solution, and what is the blow-up profile in this case. Numerical studies of blow-up in NLS equations allowed to 
give indications for  analytical results in the field, see 
\cite{SS99} for a detailed discussion and for references. 
First numerical studies of the DS II systems were presented 
in \cite{WW}, the first study of blow-ups in DS II  appeared in 
\cite{BMS}. In \cite{MFP} numerical experiments indicated  that the lump 
soliton was unstable both against blow-up and being dispersed away 
(in dependence of whether the perturbed lump has a larger or smaller 
$L^{2}$ norm than the lump). These experiments were repeated with 
much higher resolution on parallel computers in \cite{KMR}. There it 
was shown that the Ozawa solution is also unstable against being 
dispersed away or a blow-up before the blow-up time of the 
unperturbed solution. The mechanism of the latter, however, was not 
studied in \cite{KMR}. In \cite{KS}, this question was addressed in 
a similar way as blow-up in generalized Korteweg-de Vries (KdV) \cite{KP2013} and 
Kadomtsev-Petviashvili (KP) \cite{KP2014} equations. There certain norms 
of the solution were traced during the time evolution
and matched to a dynamical rescaling 
(\ref{resc}) which allowed in \cite{KP2013,KP2014} to infer the scaling 
function $L(t)$ near blow-up. Since only serial computers were available in 
\cite{KS}, the necessary resolution could not be reached in order 
to decide whether the blow-up is according to the scaling (\ref{loglog}) or 
(\ref{ozawascal}). It could be shown however that blow-up seems to be 
possible only in the integrable focusing DS II equation (\ref{DSII}). 
In general DS II equations (a factor $\beta\in\mathbb{R}$ in front of 
the $\Phi$ in the first line of (\ref{DSII}), see \cite{GS}), no blow-up  was found in 
\cite{KS}. In \cite{KR} the so-called semiclassical limit for DS II 
equations was studied numerically. It was shown that blow-up appears 
only if the initial data are (locally near the maximum) symmetric with respect to the exchange 
of the spatial coordinates $x\mapsto y$, i.e., for instance in the 
radially symmetric case. 

In the present article we study the cases with blow-up in 
\cite{KMR} and \cite{KS} with a parallized code on GPUs. This 
allows to get the necessary resolution to identify the type of 
blow-up. Concretely we study perturbed lump solutions, a perturbed 
Ozawa solution, and Gaussian initial data. The results can be 
summarized in the following conjecture:\\
\begin{conjecture}\label{conj}
 Consider initial data $\psi_{0}\in C^{\infty}(\mathbb{R}^{2})\cap 
 L^{2}(\mathbb{R}^{2})$ for the focusing DS II equation  (\ref{DSII}) 
 with a single global 
 maximum of $|\psi_{0}|$ such that the solution to DS II has a 
 blow-up in finite time. Then the blow-up is self-similar according 
 to (\ref{resc}) with a scaling factor $L(t)$ of the form 
 (\ref{ozawascal}) and the blow-up profile given by the lump, i.e.,
 \begin{equation}
      \psi(x,y,t)=\frac{P(X,Y)}{L(t)}+\tilde{\psi}, \quad
P(X,Y) = \frac{2}{1+X^{2}+Y^{2}},\quad L(t)\sim t^{*}-t
     \label{buprofile},
 \end{equation}
where  $\tilde{\psi}$ is bounded for all $t$.
\end{conjecture}

The paper is organized as follows: In section 2 we collect some 
known analytical facts on the focusing DS II equation. In section 3 we
present the used 
numerical approaches and test them for the analytically known lump 
and for Ozawa's blow-up solution. In section 4, we study the 
perturbed lump, a perturbed Ozawa solution, and Gaussian initial 
data. We add some concluding remarks in section 5.

\section{Analytical preliminaries}
In this section we review some known analytical facts on the DS II 
system which are relevant for the numerical studies to be presented 
in this paper.

For the Cauchy problem (\ref{DSII}) based on a work by Fokas and Sung 
\cite{FS}, Sung 
\cite{Su1, Su2, Su3, Su4} has proven the following      
\begin{theorem}\label{DSSung}
Let $\psi_0\in \mathcal S(\R^2)$, the space of rapidly decreasing 
smooth functions. Then \eqref {DSII} possesses a 
unique global solution $\psi$ such that the mapping $t\mapsto 
\psi(\cdot,t)$ belongs to $C^\infty(\R, \mathcal S(\R^2))$ if 
$$||\widehat{\psi_0}||_1||\widehat{\psi_0}||_\infty<C,$$ where $C$ is an 
explicitly known  constant.
\end{theorem}
Thus Sung obtains  global well-posedness  under the assumption that $\hat \psi_0\in L^1(\R^2)\cap L^\infty(\R^2)$ and $\psi_0\in L^p(\R^2)$ for some $p\in [1,2),$ see \cite{Su4}.
Perry \cite{Pe} has generalized this to the case of initial data in $H^{1,1}(\R^2)=\lbrace f\in 
L^2(\R^2) \;\text{such that}\; \nabla f, (1+|\cdot|)f \in 
L^2(\R^2)\rbrace.$ Recently Nachman, Regev and Tataru \cite{NRT}
 improved this to initial data 
$\psi_{0}\in L^{2}(\mathbb{R}^{2})$. Note that numerical results in 
\cite{KR2,KS} indicate that the bound given in Theorem 
\ref{DSSung} is not optimal, i.e., that initial data not satisfying 
the condition can still lead to global existence in time of the 
solution. 

The focusing DS II equation (\ref{DSII}) is completely integrable 
and thus  has an infinite number of conservation laws. The first two 
conserved quantities are the squared $L^2$ norm
$$ ||\psi||_2^2 =  \int_{\R^2}|\psi(x,y,t)|^2dx dy$$
and the energy (Hamiltonian)
\begin{equation}
    E(\psi(t))=\int_{\R^2}\left[
|\partial_{x}\psi|^2-|\partial_{y}\psi|^2+ (|\psi|^2+ 
\Phi)|\psi|^2)\right] dxdy
    \label{energy}.
\end{equation}

\begin{figure}[!htb]
\includegraphics[width=0.32\hsize]{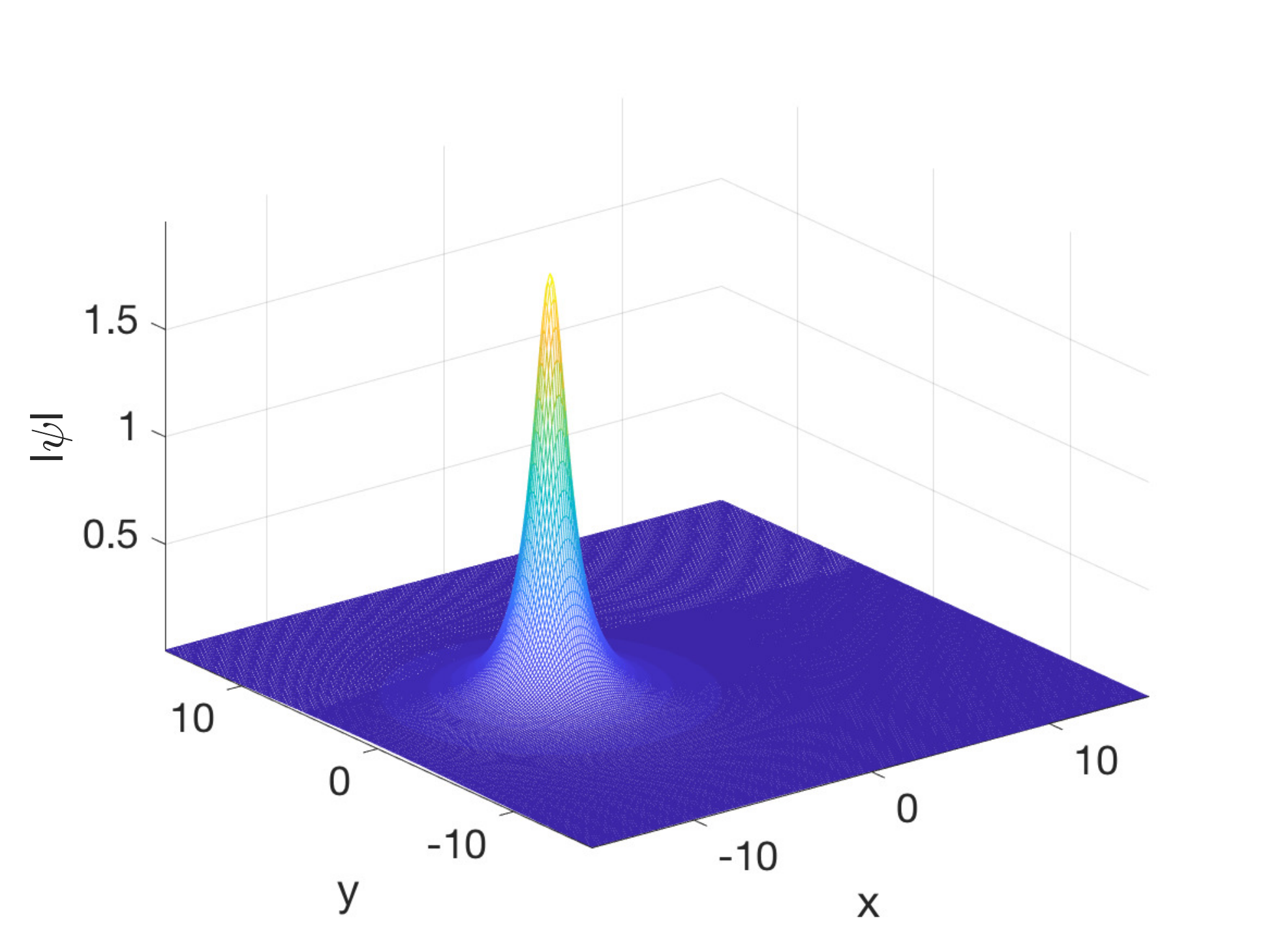}
\includegraphics[width=0.32\hsize]{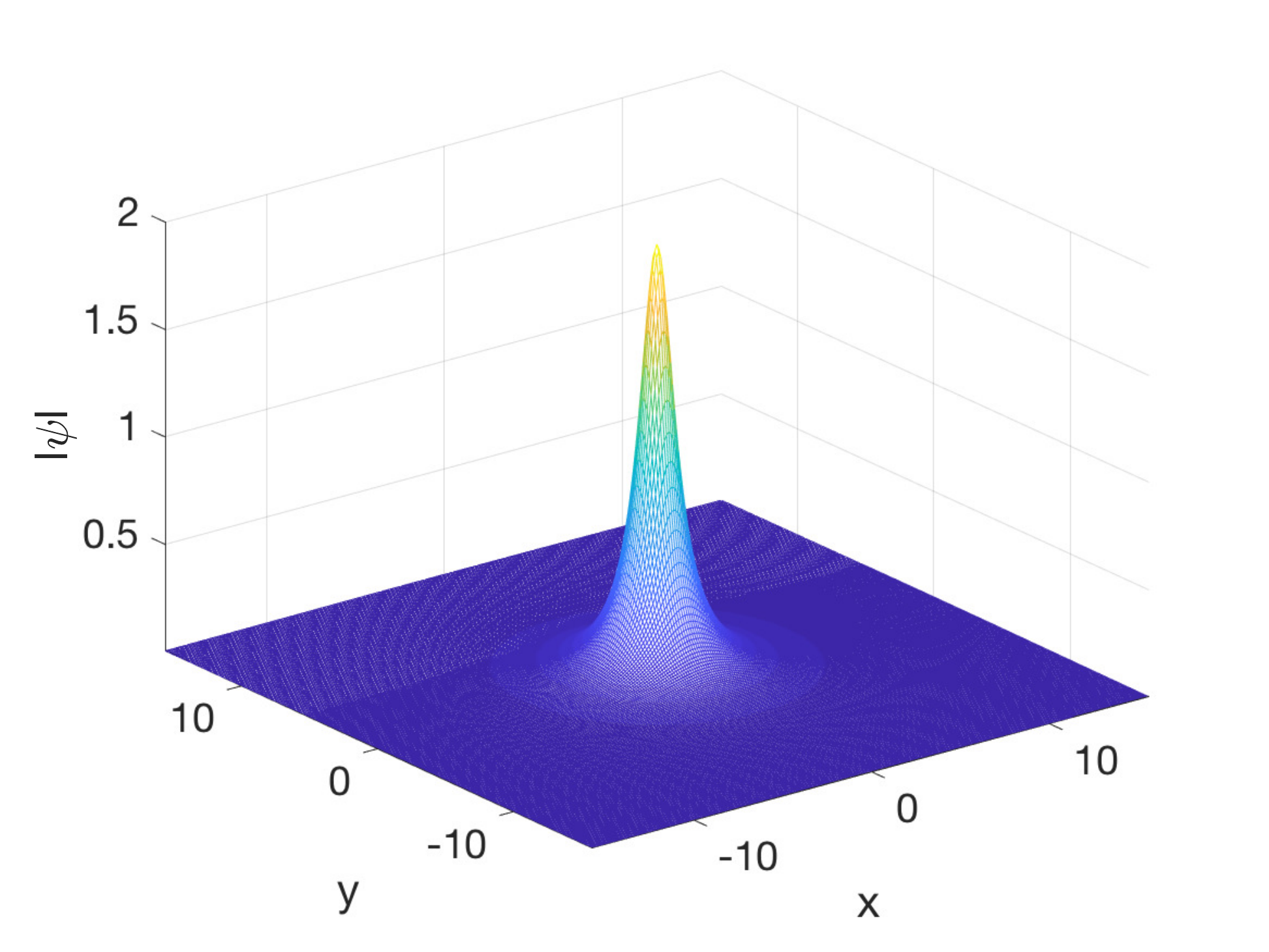}
\includegraphics[width=0.32\hsize]{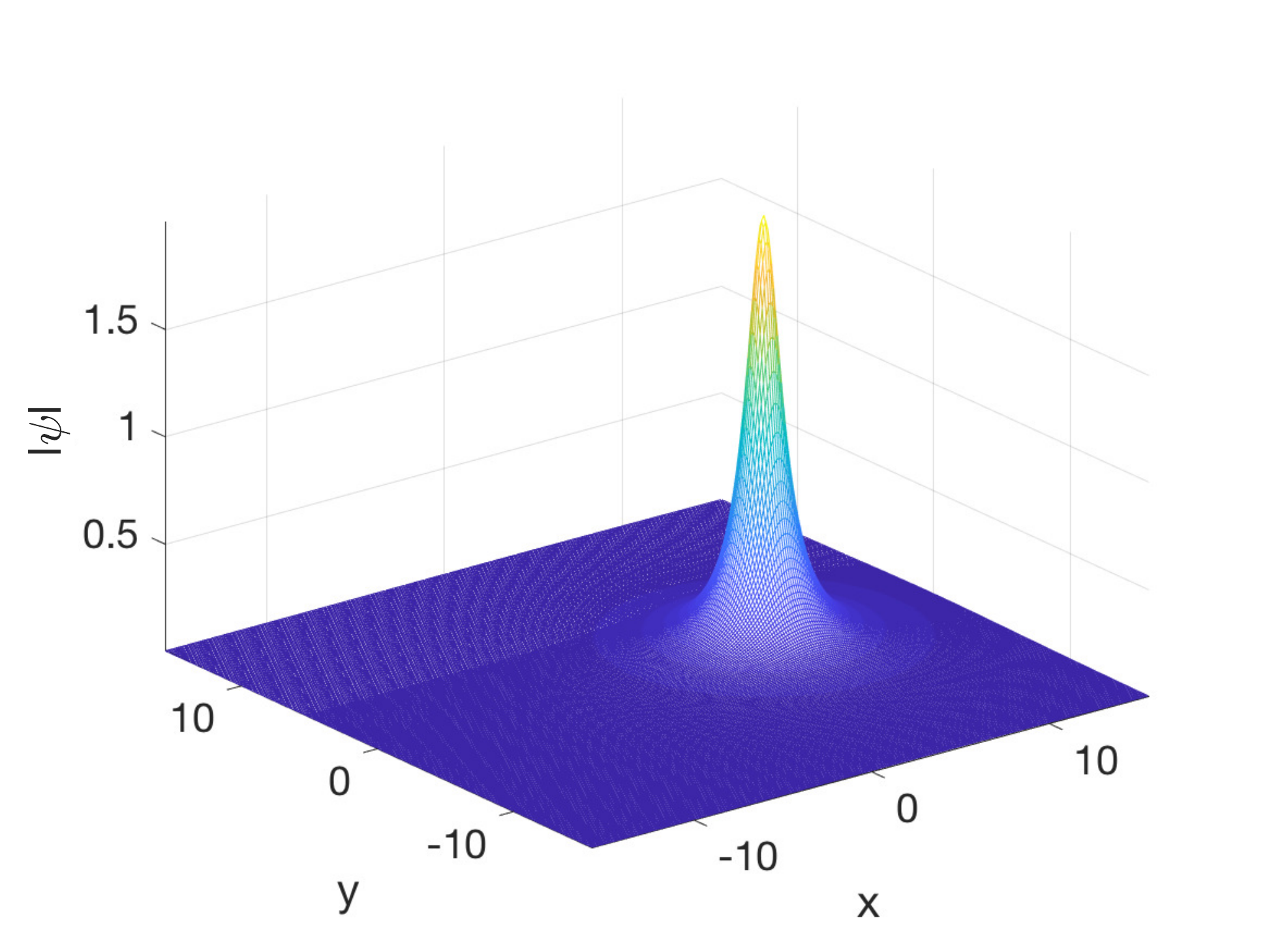}
\caption{Lump solution (\ref{lump}) with 
parameters $\xi = -1$, $\eta = 0$, $z_0 = 0$, $c = 1$: On the left 
for $t=-1$, on the middle for $t=0$, and on the right for $t=1$. }
\label{figlump3t}
\end{figure}

The focusing  DS II system possesses a family of localized solitary 
waves (\cite{APP, AC}) called lumps given by
\begin{equation}\label{lump}
 \psi(x,y,t) = 2c \frac{\exp \left( -2i(\xi x - \eta y + 
2(\xi^{2}-\eta^{2} )t)\right)}{|x + 4\xi t + i(y + 
4\eta t) +
 z_{0}|^2+|c|^2}
\end{equation}
 where $(c,z_{0})\in \mathbb{C}^2$ and $(\xi,\eta)\in\mathbb{R}^2$ 
 are constants. The lump moves with constant 
 velocity $(-4\xi, -4\eta)$ and has an algebraic decay as $(x^2+y^2)^{-1}$ for 
 $x,y\to\infty$. We show the solution for $\xi=-1$, $\eta=0$, 
 $z_{0}=0$, $c=1$ for various values of $t$ in Fig.~\ref{figlump3t}. 

Perry \cite{Pe16} showed that the lump is unstable in the inverse 
 scattering picture in the sense that a small perturbation of lump 
 initial data leads to an inverse scattering transform with different 
 analytical properties than the lump. We address perturbations of the 
 lump numerically in this paper.

\begin{figure}[!htb]
\includegraphics[width=0.32\hsize]{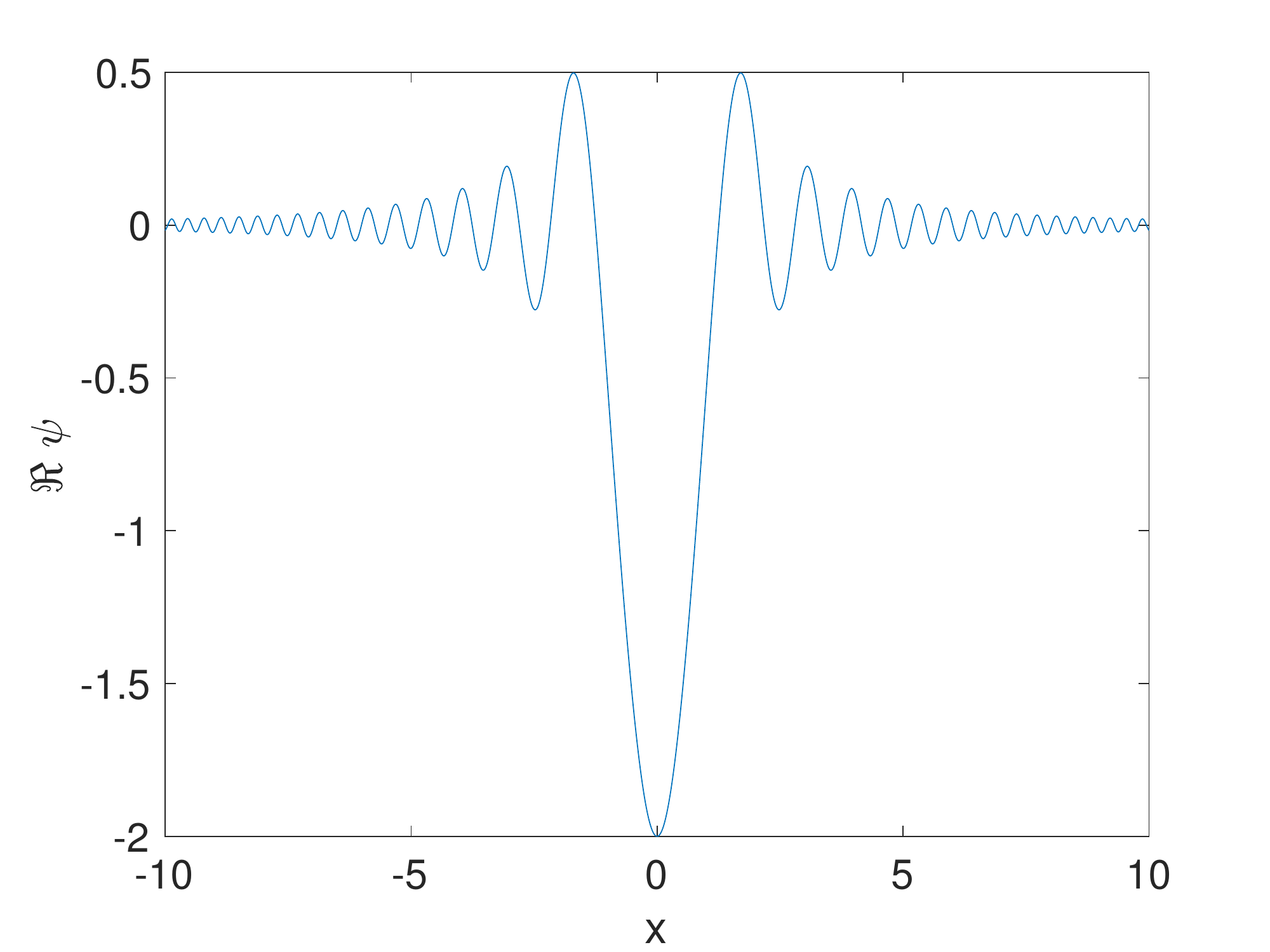}
\includegraphics[width=0.32\hsize]{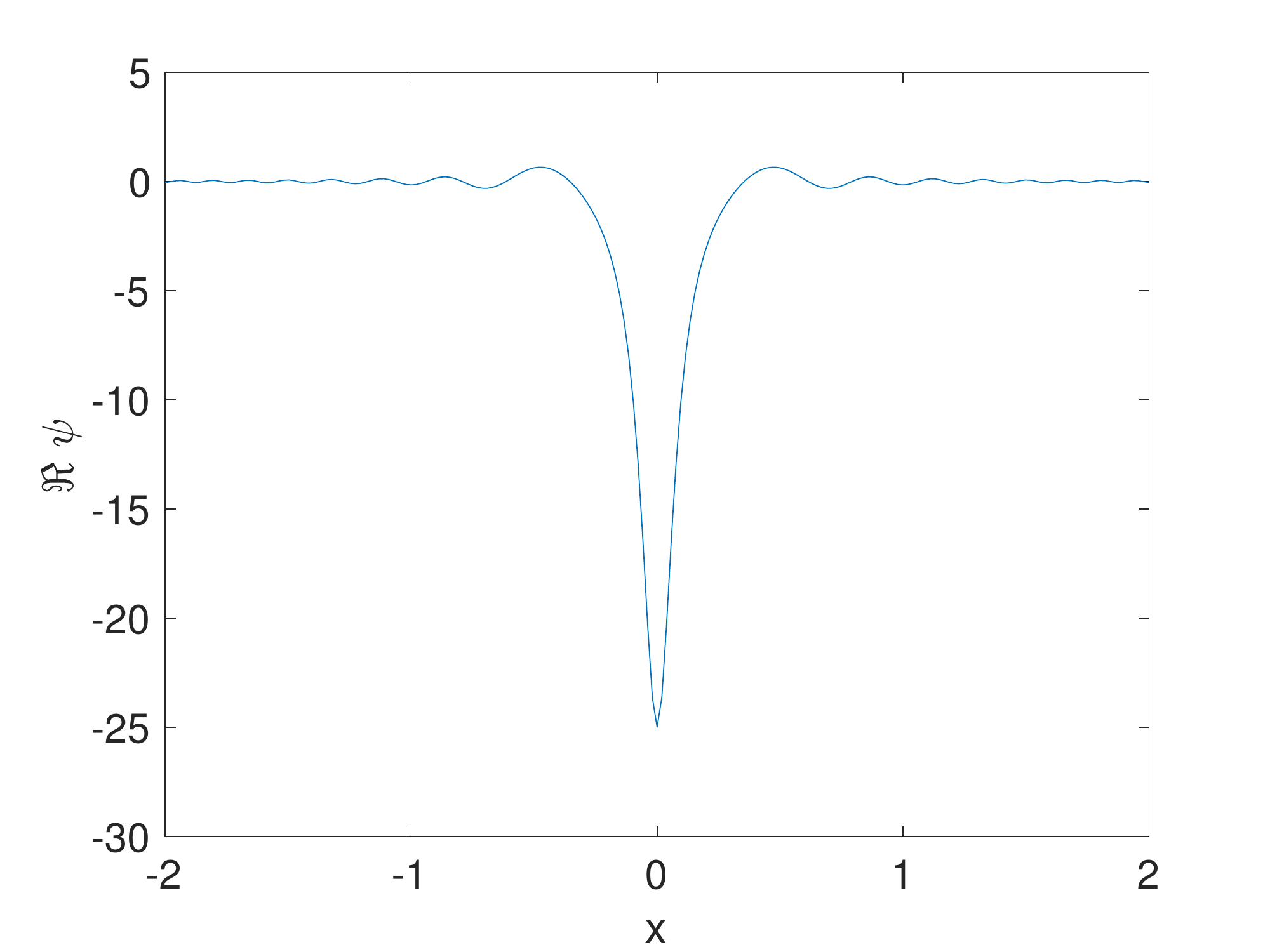}
\includegraphics[width=0.32\hsize]{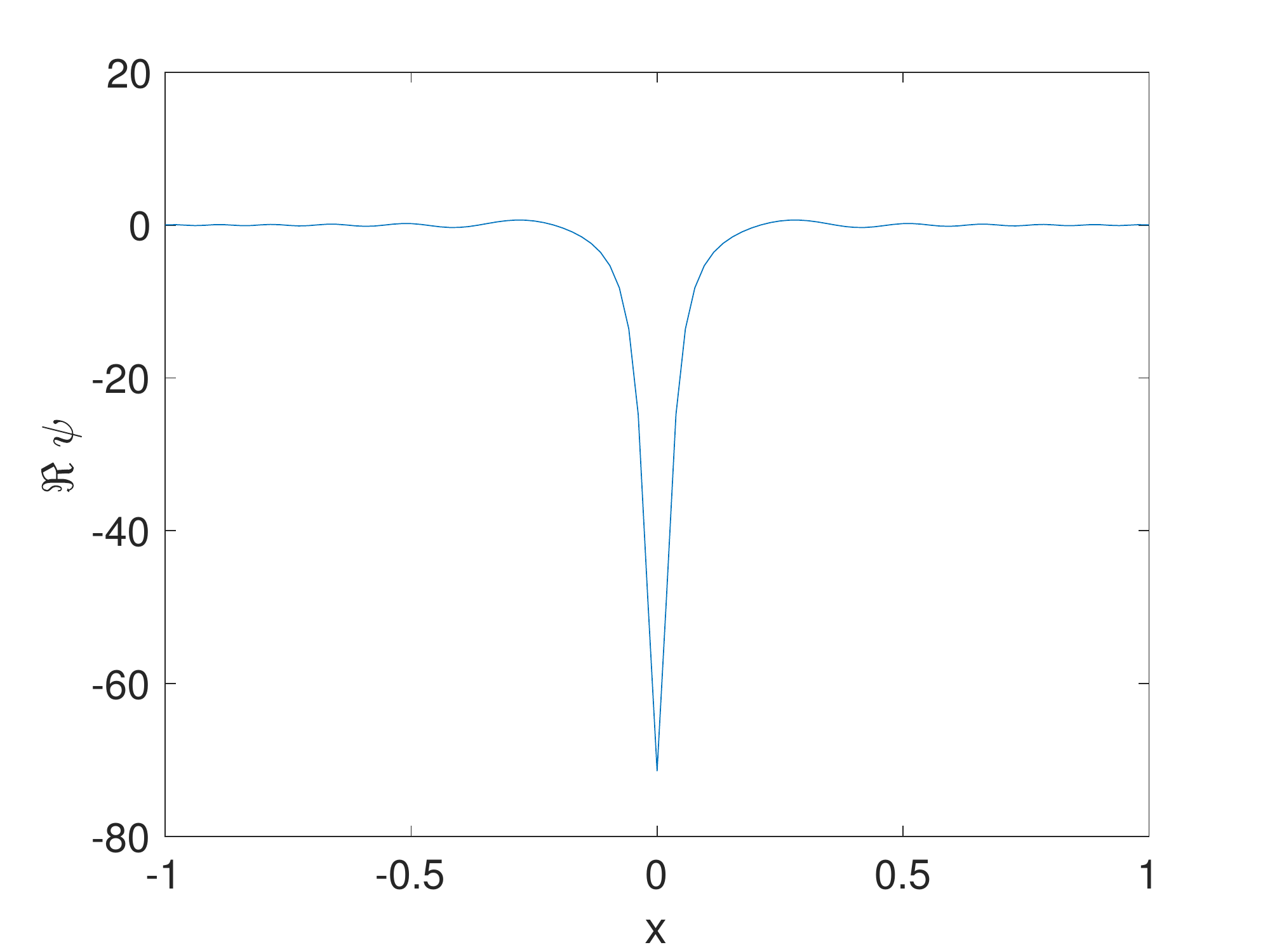}
\caption{Real part of the Ozawa solution (\ref{soluoz}) for $a = 1$, 
$b = -4$ on the $x$-axis for various values of $t$: On the left for $t = 0$, on the 
middle for $t = 0.23$ and on the right for $t = 0.243$ close to the blow-up time $t^{*} = 0.25$. }
\label{figozawa}
\end{figure}

Applying the pseudoconformal invariance (\ref{pseudo}) of DS II on the 
lump solution (\ref{lump}), Ozawa \cite{Oz} has constructed an 
explicit  blow-up solution to DS II. 
\begin{theorem}[Ozawa] 
Let   $a,b\in\R$ such that $ab<0$ and $t^{*}=-a/b$. 
Let
\begin{equation}
 \psi (x,y,t) = \exp \left( i \frac{b}{4(a+bt)} (x^2 - y^2)
\right) \frac{v(X,Y)}{a+bt}
\label{soluoz}
\end{equation}
 where 
\begin{equation*}
 v(X,Y) = \frac{2}{1+X^2+Y^2}, \,\, X=\frac{x}{a+bt},
\,\, Y=\frac{y}{a+bt}.
\end{equation*}
Then, $\psi$ is a solution of (\ref{DSII}) with 
\begin{equation*}
 ||\psi(t)||_2 = ||v||_2 = 2\sqrt{\pi}
\end{equation*}
and
\begin{equation*}
 |\psi(t)|^2 \rightarrow 4\pi \delta 
 \,\, \text{in}\quad \mathcal S'\quad \mbox{when} \,\, t \rightarrow 
 t^{*},
\end{equation*}
where $\delta$ is the Dirac measure, and where $\mathcal{S}'$ is the space 
of tempered distributions.
\end{theorem}
The mass density 
$|\psi(.,t)|^2$ of the solution converges as $t\to t^*$ to a Dirac 
measure with total  mass $||\psi(.,t)||^2_2=||\psi_0||^2_2$ . Every regularity breaks down 
at  the blow-up point, but the solution continues to solve DS II for 
$t>t^{*}$. Note that the Ozawa solution is not in $H^{1}$. 
We show the real part of the 
 solution for $a= 1$ and $b = -4$ on the $x$-axis for different values of $t$ 
 in Fig.~\ref{figozawa}.
 Note that the factor $\exp  \left(i \frac{b}{4(a+bt)}(x^{2} - 
 y^{2})\right)  = \exp  \left(-i(x^{2} - y^{2})/(t^{*} - t)\right)$ in
(\ref{soluoz}) reads $\exp \left(i(X^{2}- Y^{2} )/\tau\right)$ (for a proper choice of the integration
constants). This means that the rapid oscillations will be suppressed near
blow-up which is reached for $\tau\to\infty$. This behavior can be 
seen in Fig.~\ref{figozawa}, where the oscillations are most visible 
for $t = 0$. Note the zoom in near blow-up in $x$ as well as the 
increasing amplitude near blow-up.

\section{Numerical methods}\label{sec:num}

In this section, we review the numerical methods for the study of 
blow-up in solutions to the focusing DS II equation (\ref{DSII}). The 
approaches are tested via the propagation of lumps and  the Ozawa 
solution. 

\subsection{Numerical methods for the time-evolution}\label{sec:numtime}

In order to  numerically integrate equation \eqref{DSII}, we use
a Fourier spectral method in the spatial coordinates. This means 
that the standard Fourier transform is approximated via a discrete 
Fourier transform which can be efficiently computed via the 
\emph{fast Fourier transform} (FFT). In an abuse of notation we will 
use the same symbols for the continuous and the discrete Fourier 
transform.
The discretization in Fourier space implies that (\ref{box}) is 
approximated via a finite dimensional system of 
 ordinary differential equations for the Fourier 
coefficients $\hat{\psi}=\mathcal{F}\psi$ of $\psi$ of the form
\begin{equation}
   \partial_t \widehat{\psi}=\mathcal{L}\widehat{\psi}+\mathcal{N}(\psi),
    \label{fnlsfourier}
\end{equation}
where $\mathcal{L}=-i(\xi_{1}^{2}-\xi_{2}^{2})$, and where
$\mathcal{N}(\psi)=2i \mathcal{F}[\mathcal{F}^{-1}\{(\xi_{1}^{2}-\xi_{2}^{2})/(\xi_{1}^{2}+\xi_{2}^{2})
 \mathcal{F}(|\psi|^2)\}\psi]$ 
denotes the nonlinearity; here $\xi_{1}$, $\xi_{2}$ are the dual 
Fourier variables for $x$ and $y$ respectively.  For systems of 
the form (\ref{fnlsfourier}) with diagonal $\mathcal{L}$ as in the 
case of Fourier methods, there are many efficient high-order time 
integrators, see \cite{KR} for a comparison of fourth order methods for  DS II.

As in \cite{KMR}, we use here a fourth order splitting method.
The  motivation for these methods is the Trotter-Kato formula
\cite{TK,Ka}
\begin{equation*}
 {\lim}_{n\rightarrow\infty}\left(e^{-tA/n}e^{-tB/n}\right)^{n}=e^{-t\left(A+B\right)},
\end{equation*}
where $A$ and $B$ are certain unbounded linear operators, for details 
see \cite{Ka}. Splitting techniques have been studied by Bagrinovskii and
Godunov in \cite{BG}, by Strang \cite{ST}, and for hyperbolic 
equations in Tappert \cite{Tap} and Hardin and Tappert \cite{HT}
(split step method for NLS).

For an equation of the form $u_{t}=\left(A+B\right)u$ the idea is to write the solution in the
form
\[
u(t)=\exp(c_{1}tA)\exp(d_{1}tB)\exp(c_{2}tA)\exp(d_{2}tB)\cdots\exp(c_{k}tA)\exp(d_{k}tB)u(0),
\]
where the real numbers $(c_{1},\,\ldots,\, c_{k})$ and $(d_{1},\,\ldots,\, d_{k})$
represent fractional time steps.  
Yoshida \cite{Y} gave an approach 
which produces split step methods of any even order. 
Here the DS II equation (\ref{fnlsfourier}) is split into
\begin{eqnarray}\label{2.4a}
\partial_{t} \hat{\psi}=-i(\xi_{1}^{2}-\xi_{2}^{2})\hat{\psi} ,
\\
\partial_{t} \hat{\psi}= 
2i\mathcal{F}\left[\mathcal{F}^{-1}\left\{\frac{\xi_{1}^{2}-\xi_{2}^{2}}{\xi_{1}^{2}+\xi_{2}^{2}}\mathcal{F}\left|\psi\right|^{2}\right\}\psi\right],  
\label{2.4}              
\end{eqnarray}
which are both explicitly integrable, the first one in Fourier space, 
the second in physical space since 
$|\psi|^{2}$  is constant in time for (\ref{2.4}). Because of this 
property and the linearity of equation (\ref{2.4a}), the splitting 
scheme conserves the $L^{2}$ norm of $\psi$ which implies that the 
method is stable. But the conservation of the $L^{2}$ norm by the 
splitting scheme also has the consequence that the latter cannot be used to control 
the numerical accuracy. Instead we consider 
the numerically computed energy (\ref{energy}).
It was shown in \cite{KR} that the quantity 
$\Delta_{E}=E(t)/E(0)-1$ can be used in this case to indicate
the numerical accuracy which is overestimated  by one to  two 
orders of magnitude.

\subsection{Singularity tracing and self-similar blow-up}
It is well known that  a (single) singularity $z_0\in \C$ of a 
function $f$ of the form $f(z) \sim (z-z_{0})^{\mu}$, with $\mu\not 
\in \mathbb Z$ implies  the following asymptotic behavior for the corresponding Fourier transform
\begin{equation}
    |\widehat{f}(k)|\sim 
    \frac{1}{k^{\mu+1}} e^{-k\delta},\quad |k| \gg 1,
    \label{fourasymp}
\end{equation}
where $\delta=\text{Im}\, z_{0}$. 
If numerically the Fourier transform is approached via an FFT, 
relation (\ref{fourasymp}) can still be used to identify the type of 
the singularity via the quantity $\mu$ obtained by a fitting of the 
Fourier coefficients for large wave numbers to  
(\ref{fourasymp}). This was first used in \cite{SSF} to identify the 
appearance of a singularity. 

In  \cite{KR2013a,KR2} this approach was discussed with considerably higher resolution in detail and 
 used to quantitatively identify the time where 
the distance $\delta$ of the singularity from the real axis becomes 
smaller than the minimal resolved 
distance via Fourier methods, i.e., $ m:=2\pi D/N$
with $N\in \N$ being the number of 
 Fourier modes and $2\pi D$ the length of the computational domain in 
 physical space.  A value of $\delta <m$ 
 cannot be distinguished numerically from $0$.

Unfortunately this method is not applicable here for most of the 
studied examples because of the slow algebraic decrease of the lump 
and the Ozawa solution towards infinity. This implies that the 
approximation as a periodic function is as discussed 
not analytic on the domain boundaries. This leads already to an 
algebraic decrease of the Fourier coefficients, an effect which is 
not easily separated from the to be identified singularity in the 
complex plane. The Ozawa solution has the additional problem that it 
is rapidly oscillating which as discussed below leads also to 
oscillations in the Fourier coefficients. In addition the singularity 
tracing requires a very high resolution in time which is difficult to 
achieve in two-dimensional computations. Thus we proceed here as in 
\cite{KMR} and stop the code once the computed energy  indicates 
that the accuracy drops below plotting accuracy (typically if the 
quantity $\Delta_{E}<10^{-3}$). 
 
The  time at which the code is stopped because of the above criterion 
is obviously not the blow-up time itself.  The latter will be 
identified by tracing the $L^{\infty}$ norm of the solution and 
matching it 
to what would be expected from a dynamical rescaling (\ref{resc}) of 
the DS II equation. Assuming that $L(t)$ vanishes at the blow-up, we 
can study the type of the 
blow-up for DS II in a similar way as it has been done for 
generalized KdV equations in \cite{KP2013}, for generalized KP equations in   
\cite{KP2014} and for fractional NLS equations in \cite{KSM}: we 
integrate DS II directly, as described above, and then we use some 
post-processing to characterize the type of blow-up via the  
rescaling (\ref{resc}). Since the $L^{2}$ norm of $\psi$ is invariant 
under this rescaling, we consider the $L^{\infty}$ norm of $\psi$. We get with (\ref{resc}) 
that 
\begin{equation}
    ||\psi(t,\cdot) ||_{\infty}\propto \frac{1}{L}   
    \label{genscal}
\end{equation}
for $t\to t^{*}$. 
Thus tracing the norm and matching it to a factor 
$(t^{*}-t)^{\gamma}$, it should be possible to determine both the 
blow-up time $t^{*}$ and the  exponent $\gamma$. Concretely we trace 
the $L^{\infty}$ norm during the time evolution until the results as 
indicated by the energy conservation can no longer be trusted. 
Then we fit $||\psi||_{\infty}$, according to
\begin{equation}
    \ln ||\psi||_{\infty}\sim \alpha+ \gamma \ln (t^{*}-t),
    \label{logfit}
\end{equation}
for a certain number of the last recorded time steps before the code 
is stopped. The fitting is done with the optimization algorithm 
\cite{fminsearch} distributed with Matlab to find constants $\alpha$, $\gamma$ 
and $t^{*}$ such that the residual 
of $\ln ||\psi||_{\infty}-[\alpha+ \gamma \ln (t^{*}-t)]$ becomes 
minimal. 
Note that it will be numerically impossible to distinguish the case 
$L(t)\propto \sqrt{t^{*}-t}$ from the formula (\ref{loglog}) established analytically 
for $L^{2}$ critical NLS equations since the effect of the logarithms will 
be too small. Note also that the $L^{2}$ norm of $\psi_{x}$, which 
has been applied in \cite{KP2013} and \cite{KP2014} in a similar 
context,  cannot be used 
in the case of the Ozawa solution for which the latter is not finite. 
Though here we approximate this solution  via a periodic function for 
which all these norms are finite, the divergence of 
$||\psi_{x}||_{2}$ on the whole real line leads to numerical problems 
on the considered large computational domains. 
Therefore we only consider $||\psi||_{\infty}$ here.

The rescaling (\ref{resc}) implies for the equation (\ref{box})
\begin{equation}
    \begin{split}
    i\partial_{\tau}\Psi+(\partial_{X}^{2}-\partial_{Y}^{2})\Psi
    +i\partial_{\tau}(\ln 
    L)(X\partial_{X}\Psi+Y\partial_{Y}\Psi+\Psi)+&\\
    2\left[(\partial_{X}^{2}+\partial_{Y}^{2})^{-1}(\partial_{X}^{2}-\partial_{Y}^{2}))
    |\Psi|^{2}\right]\Psi&=0.
    \end{split}
    \label{boxresc}
\end{equation}
Whereas it is numerically inconvenient to solve this equation instead 
of DS for our examples (the lump and the Ozawa solution have 
algebraic decay towards infinity), it is interesting for theoretical 
reasons. In the limit $\tau\to\infty$, the blow-up is reached, and it 
is assumed that $\Psi$ is $\tau$ independent in the limit. If $L 
\propto t^{*}-t$ as for 
the Ozawa solution, also $\partial_{\tau}\ln L$ would vanish which 
leads for (\ref{boxresc}) to   
\begin{equation}
   (\partial_{X}^{2}-\partial_{Y}^{2})\Psi+
    2\left[(\partial_{X}^{2}+\partial_{Y}^{2})^{-1}(\partial_{X}^{2}-\partial_{Y}^{2}))
    |\Psi|^{2}\right]\Psi=0.
    \label{boxasym}
\end{equation}
The only known stationary, localized solution to this equation
appears to be the lump though it is unclear whether it is the unique solution
with this property. This would imply that a dynamically rescaled (according
to (\ref{resc})) lump could provide the blow-up profile for focusing DS II solutions.

\subsection{Parallel computing on Graphic Processor Units }
Since in the context of solutions of the focusing DS II equation 
strong gradients and rapid oscillations have to be resolved in two 
dimensions, the required high resolution cannot be reached by 
single thread codes. For this reason we use parallel computing in the 
form of heterogeneous computing with Graphics Processing Unit (GPUs). The 
code is implemented via NVIDIA's CUDA C/C++ platform. Using CPU multithreading the 
parallelized stepper is run on the GPUs, whereas the control 
parameters  such as the energy, $L^2$ norm etc.\ are concurrently 
computed on the CPU using the standard FFTW\footnote{fftw.org} 
library, offloading computational duties from the GPUs. 

The codes are run on a desktop computer equipped with two 
NVIDIA Quadro M5000s. The maximum resolution possible with the described hardware and all data used by the solver kept on the GPUs  is $2^{14} \times 
2^{14}$. All examples are run under the same scheme - first we do 
2000 time steps to get close to the blow up point, followed by 2000 
steps of size $\mathbb{d}t = 10^{-5}$ (except in one case, as described below), that run through the blow up 
point. As already explained the code cannot be trusted after the 
distance between the pole in the complex plane and the real line 
becomes smaller than the spatial resolution or the conservation of 
the numerically computed energy drops below an acceptable threshold.                       
Thus we perform exploratory low resolutions runs to establish 
roughly the blow-up time, and then apply the above approach for a high 
resolution analysis.

\subsection{Tests}\label{tests}
We test the numerical code on two examples. First, we propagate the 
analytically known 
lump solution, and then we simulate the Ozawa solution in order to 
test our approaches to identify blow-up. 

\subsubsection{Lump solution}
\begin{figure}[!htb]
\includegraphics[width=0.45\hsize]{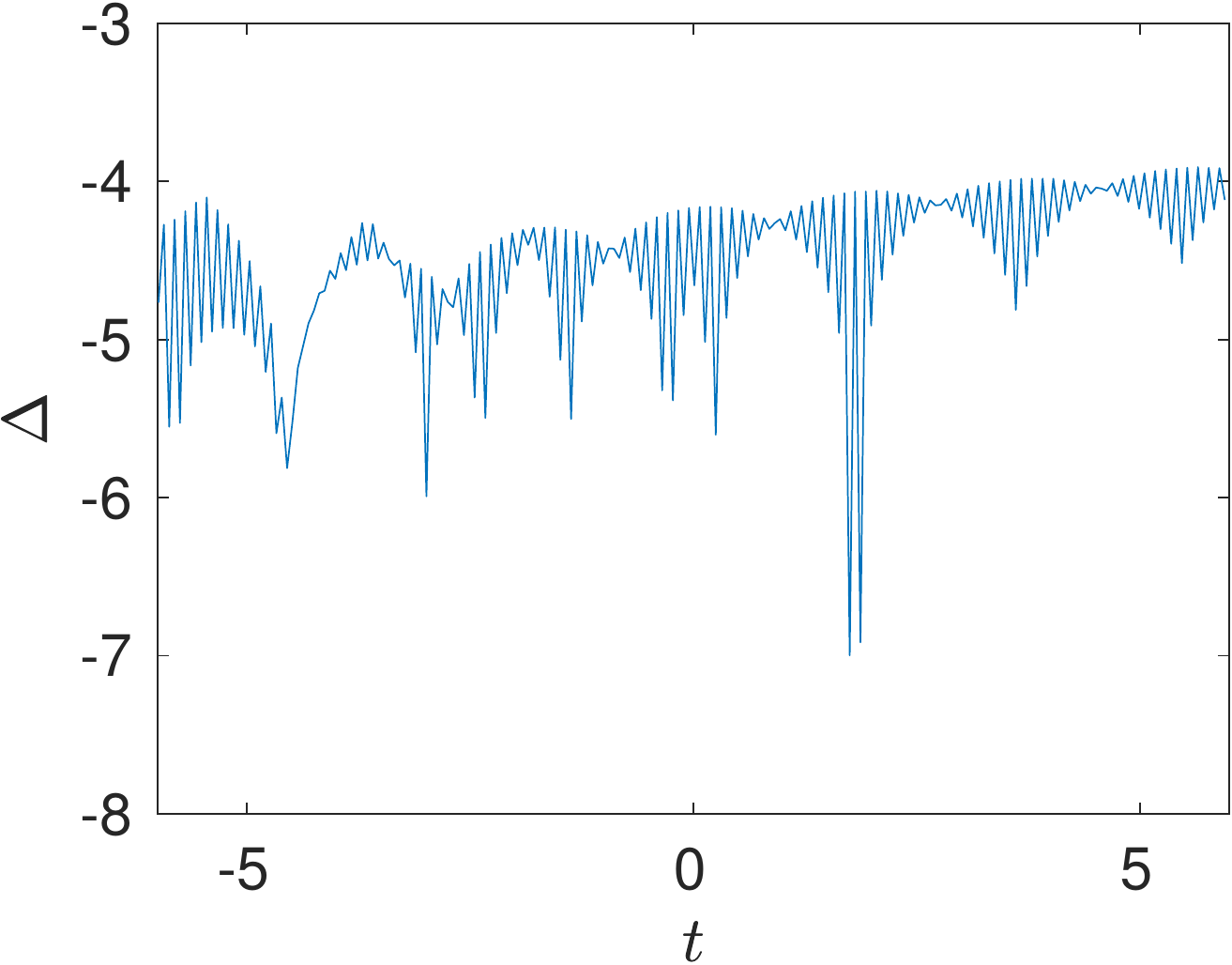}
\includegraphics[width=0.45\hsize]{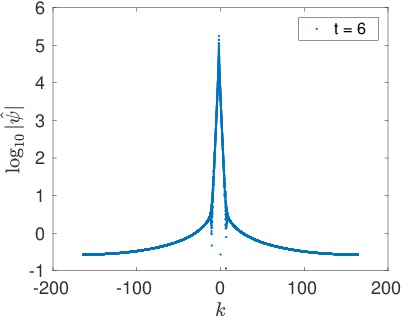}
\caption{Propagation of the lump solution (\ref{lump}) with 
parameters $\xi = 0$, $\eta = -1$, $z_0 = 0$, $c = 1$: On the left 
the quantity 
$\Delta:=||\psi_{num}-\psi_{exact}||_{\infty}/||\psi_{exact}||_{\infty}$
indicating the numerical error, on the right the modulus of the Fourier 
coefficients on the $\xi_{1}$-axis  at $t=6$. For the energy preservation we have 
$\max\Delta E = 10^{-9.2994}$. }

\label{figlump1}
\end{figure}

We first test the lump solution (\ref{lump}). We use a domain with $D = 50$ (we remind that the domain size is 
$2\pi D \times 2\pi D$) in order to handle the algebraic fall off and 
solution parameters $\xi = 0$, $\eta = -1$, $z_0 = 0$, $c = 1$, and
$t \in [-6,6]$, i.e., 
\begin{equation*}
\psi(x,y,-6) = 2\frac{\exp(-2i(y + 12))}{|x + i(y + 24)|^2 + 1}
\end{equation*}
and propagate the solution to $t = 6$ in $N_t=1000$ temporal steps, 
comparing it to the exact solution. The Fourier coefficients for the 
solution (at the final time, but the difference to the coefficients 
at the initial time is hardly visible) are shown on the right of 
Fig.~\ref{figlump1}. Due to the algebraic decay of the solution 
towards infinity, the solution is not analytical at the domain 
boundaries as a periodic function. This leads to the algebraic 
decay of the Fourier coefficients with the wave number $k$ for $|k|$ 
to infinity. The choice of the parameter $D$ in the computational 
domain $D[-\pi,\pi]\times D[-\pi,\pi]$ is motivated by exactly this 
decrease of the Fourier coefficients: for small values of $D$, the 
resolution is high in physical space, but the lack of analyticity 
leads to a slower decrease of the coefficients and thus of the 
numerical error. The latter problem is reduced for larger values of 
$D$, but then the number of points in the FFT per interval becomes 
smaller leading again to larger numerical errors. The used value 
$D=50$ is a compromise between these two effects. 

The relative decrease of the Fourier coefficients by roughly 5 orders 
of magnitude (see Fig.~\ref{figlump1} on the right) indicates that a numerical error of the order of 
$10^{-5}$ can be expected in this case. That this is indeed the case 
can be seen on the left of Fig.~\ref{figlump1} where 
the $L^{\infty}$ norm of the difference between numerical and exact 
solution (normalized by the $L^{\infty}$ norm of the exact solution) 
can be seen in dependence of time. The computed relative energy is 
conserved to the order of $10^{-10}$. This shows that the problem is 
much better resolved in time than in space, the reason being as 
discussed the algebraic decay of the solution towards infinity.


\subsubsection{Ozawa solution}
\begin{figure}[!htb]
\includegraphics[width=0.45\hsize]{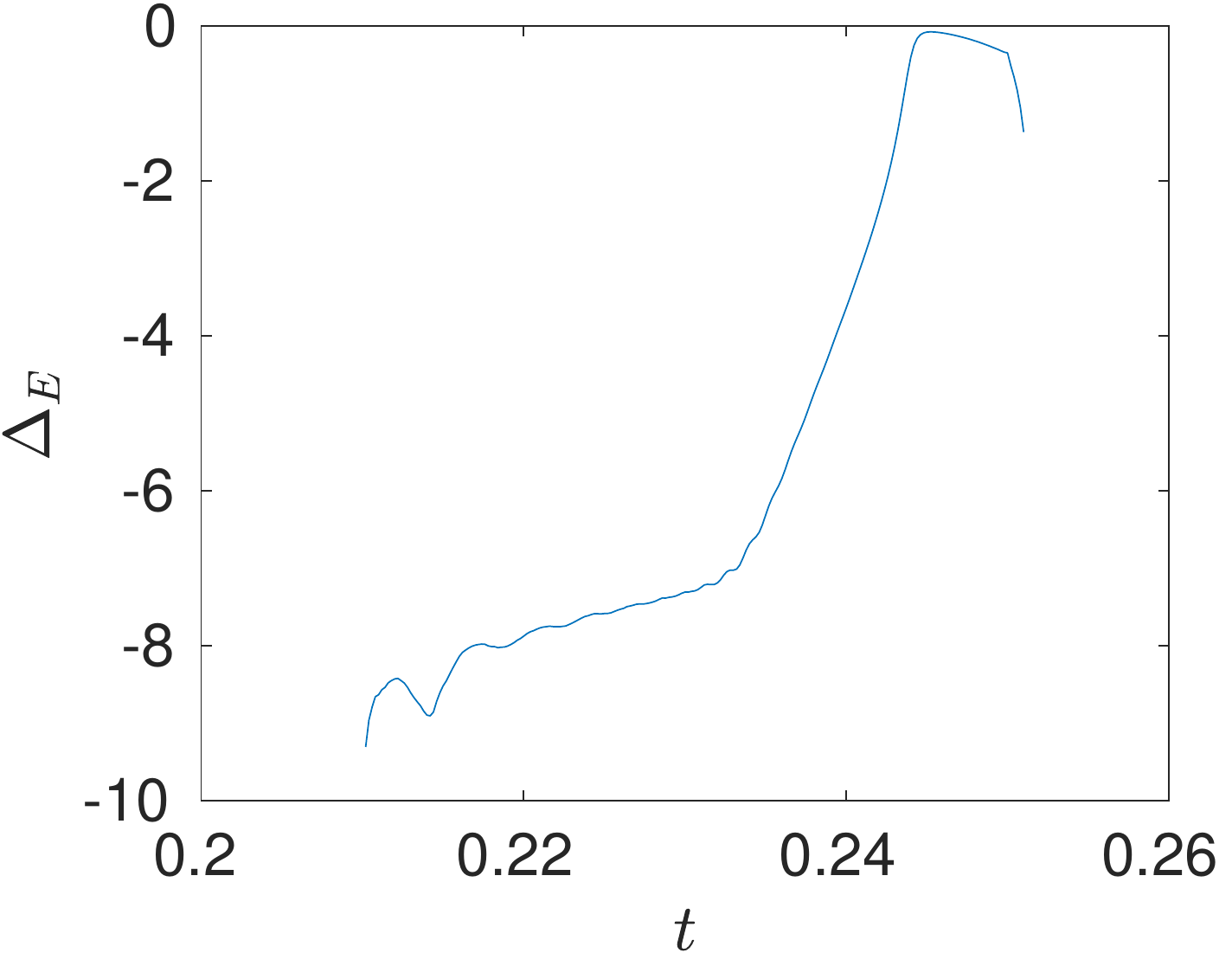}
\includegraphics[width=0.45\hsize]{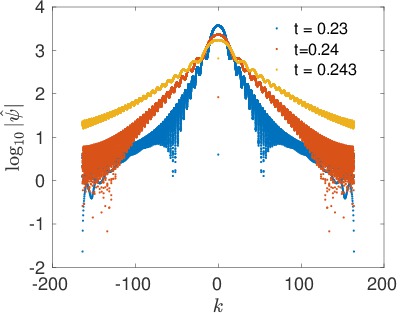}
\caption{Modulus of the Fourier coefficients for the Ozawa solution 
with parameters $a = 1$, $b = -4$ on the $\xi_{1}$-axis, for various values of $t$ on the 
right, and  
the corresponding relative computed energy $\Delta E$ on the left. }
\label{ozfourier}
\end{figure}

We consider the Ozawa solution (\ref{soluoz}) on a domain with $D = 50$ with parameters $ a = 1$ and $b = -4$ having a blow-up at $t^* = 0.25$
i.e., the situation shown in Fig.~\ref{figozawa}. For $t = 0$ this leads to the initial condition 
$$\psi(x,y,0)=2\frac{\exp(-i(x^{2}-y^{2}))}{1+x^{2}+y^{2}}.$$
The algebraic decay of the solution is clearly seen in the Fourier coefficients in
Fig.~\ref{ozfourier}, on the right for different times. The 
highly oscillatory character of the solution together with the 
algebraic decay already present for the lump makes these initial data 
especially challenging for the numerics. It can be seen that the 
oscillations in physical space lead together with the slow decay of the 
solution towards infinity to oscillations in coefficient space. Since 
the oscillations near blow-up are less important as discussed in the 
previous section, also the oscillations in coefficient space are 
reduced at later times.   

It is well known (see already \cite{BMS} and \cite{KMR}) that splitting codes for DS
will continue to run even after an analytical blow-up is encountered. This is
different from generalized KdV codes where in general overflow errors lead to
a stop of the code. However, some quantities show a jump at the numerically
encountered blow-up. Since the split-step method preserves the $L^{2}$ norm
exactly, the latter is essentially continuous there, but we see a jump in the
energy in Fig.~\ref{ozfourier} on the left. We start our simulation at $t=0.21$ and execute 4000 steps with step 
size $\mathbb{d}t = 10^{-5}$. The fitting of the  Fourier coefficients 
shows that the critical distance is reached at $t = 0.2430$ (step 
3330). However, the energy preservation  reaches the critical point 
$\Delta_E = 10^{-4}$ before that, $t = 0.2380$ at step 2980. 

\begin{figure}[!htb]
\includegraphics[width=0.45\hsize]{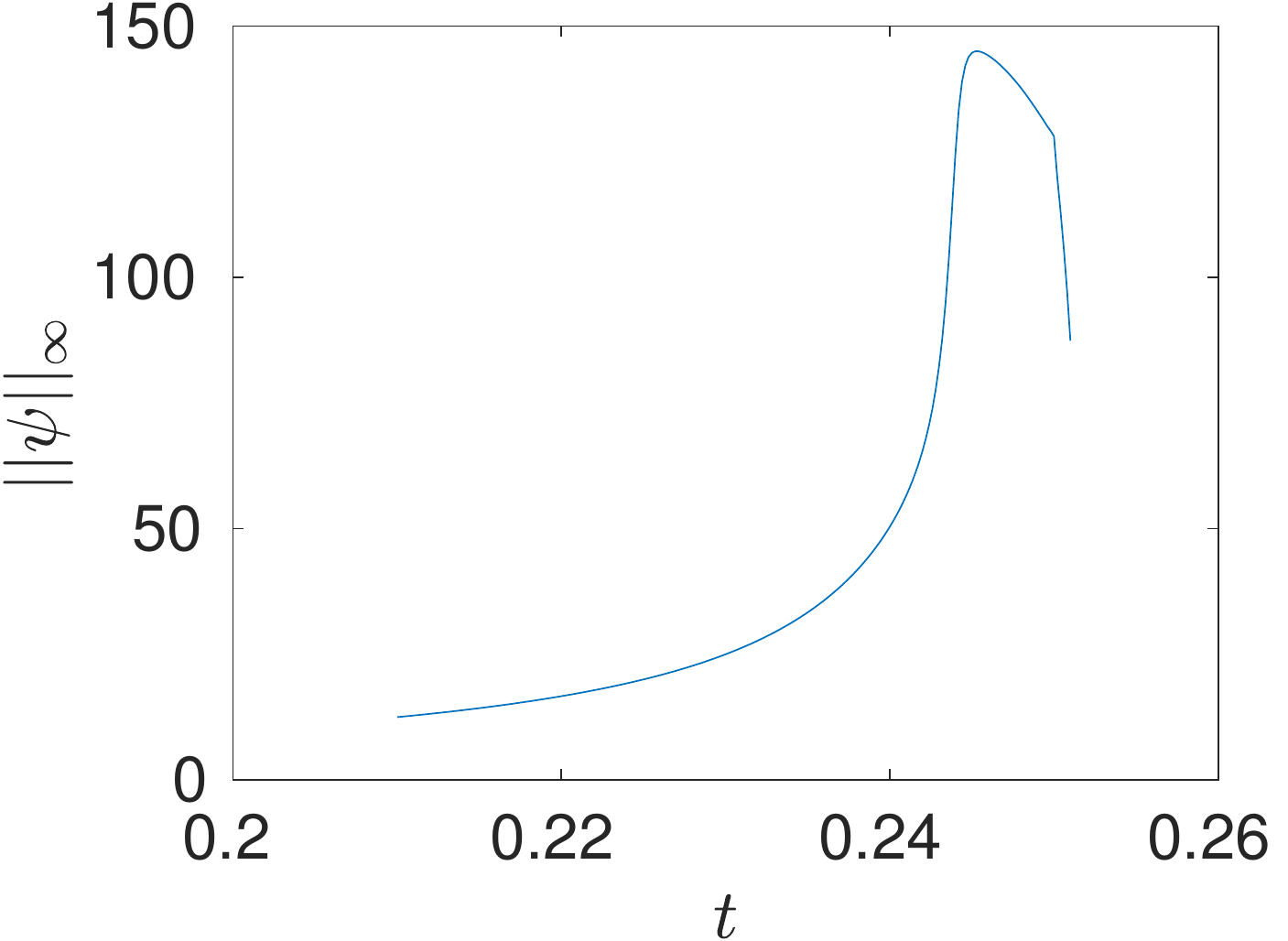}
\includegraphics[width=0.45\hsize]{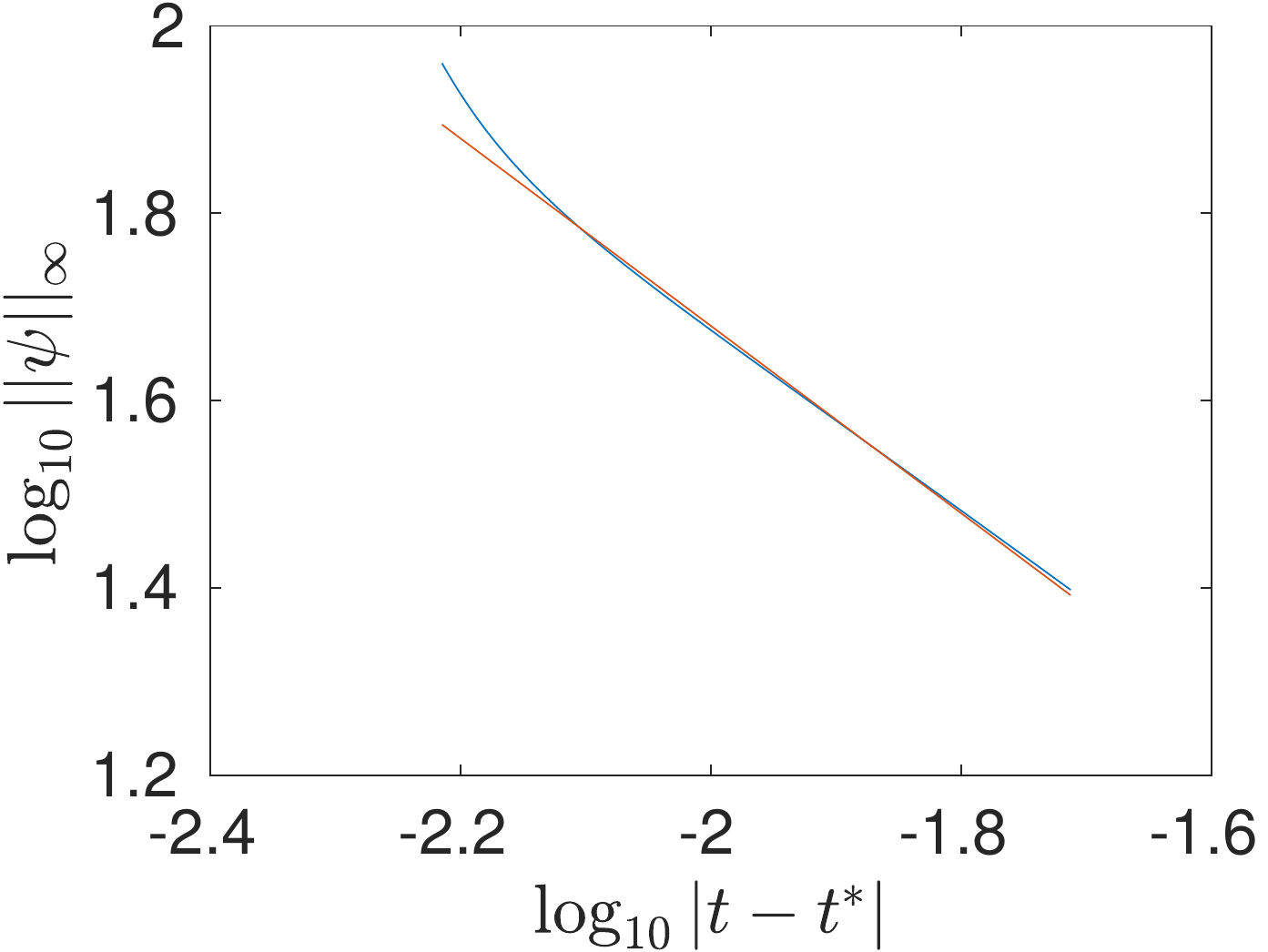}
\caption{Numerically computed Ozawa solution with parameters $a=1$, 
$b=-4$, with blow-up time $t^{*} =0.25$ on a domain with $D = 50$: on 
the left
$||\psi||_{\infty}$ in dependence of time,  on the right 
the $L^{\infty}$ norm  
 for the last 1000 recorded time steps 
in a log-log plot together with a line with  slope $-1$.}
\label{figOz}
\end{figure}

In Fig.~\ref{figOz} we show the $L^{\infty}$ 
norm of the solution on the left. Fitting this norm for
the last 1000 time steps before losing resolution to formula 
(\ref{logfit}), we get the blow-up 
time $t^{*}= 
0.2494$ and the blow-up rate of $\gamma = -0.9647$ as can be seen on 
the right of Fig.~\ref{figOz}, both in excellent agreement with the 
analytic values $0.25$ and $-1$ respectively.     

\section{Numerical results}

In this section we study blow-up in solutions to the focusing DS II 
equation (\ref{DSII}) for various examples:  two different instances  of 
perturbed  Ozawa initial data, perturbed lump initial data and 
Gaussian initial data in the so-called semi-classical regime. All 
examples except the Gaussian are computed on a domain with $D = 50$, 
due to their algebraic fall off at infinity. For the Gaussian we use 
$D = 2$, as the exponential fall off allows the solution to be well 
resolved on a much smaller domain.  The results of this section are 
summarized in Conjecture \ref{conj}.

\subsection{Ozawa Solution deformed by a small Gaussian}
\begin{figure}[!htb]
\includegraphics[width=0.45\hsize]{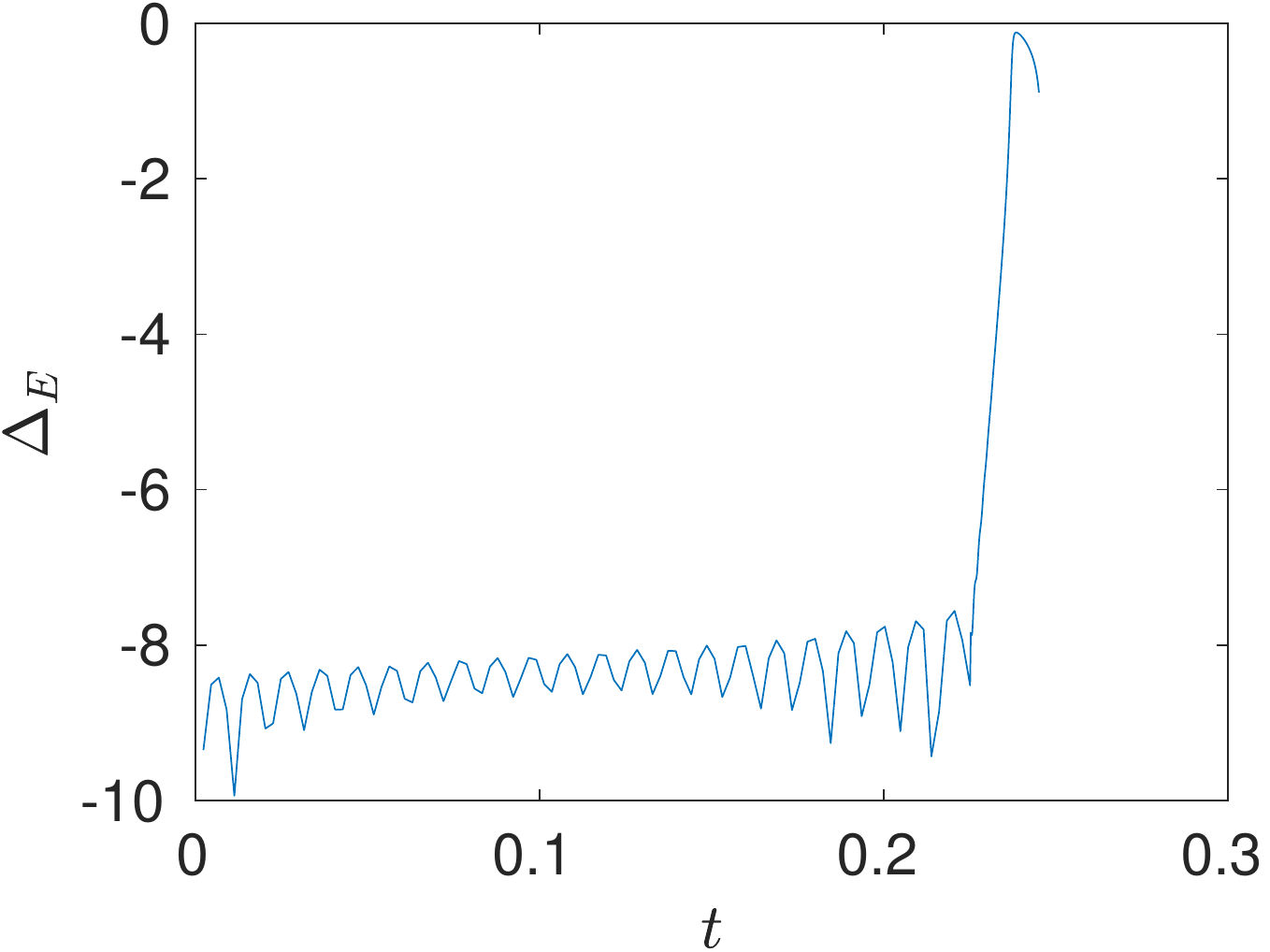}
\includegraphics[width=0.45\hsize]{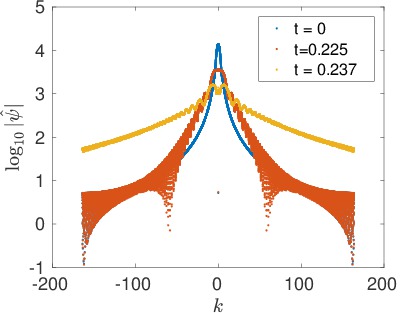}
\caption{Solution to the DS II equation  for Ozawa initial data  
deformed by a small Gaussian $0.1\exp(-x^2 -y^2)$: on the left the relative computed energy $\Delta E$, on 
the right the modulus of the Fourier coefficients  on the 
$\xi_{1}$-axis at three 
different times. }
\label{figOzdef1}
\end{figure}

We consider Ozawa initial data with parameters $ a = 1$ and $b = -4$, and deform it by 
adding a Gaussian centered at the origin of the coordinate 
system\footnote{Note that there are perturbations of the Ozawa 
solution with smaller $L^{2}$ norm of the initial data than for 
Ozawa leading to a solution without blow-up, see \cite{KMR}; here we 
are interested in the blow-up mechanism and thus only consider 
examples with blow-up.}  
\begin{equation*}
\psi(x,y,0) = 2\frac{\exp(-i(x^2 - y^2))}{1 + x^2 + y^2} + 0.1\exp{(- 
x^2 - y^2)}.
\end{equation*}
The computation is run for 2000 steps to time $0.225$ and then for 
another 2000 steps with step size $\mathbb{d}t = 10^{-5}$. The relative computed energy  $\Delta_E$ 
shown on the left of Fig.~\ref{figOzdef1}
has a jump at $t \approx 0.2326$ (step 2760) which indicates a loss 
of resolution in time.     Therefore data for $t>0.2326$ will not be considered in the following 
though the Fourier coefficients shown in Fig.~\ref{figOzdef1} on the 
right still indicate a relative spatial resolution of better than plotting 
accuracy. A fitting of the Fourier coefficients according to 
formula (\ref{fourasymp}) indicates an approaching of the singularity 
in the complex plane only at a later time. This means that we run in 
this case first out of resolution in time before losing resolution in 
the spatial coordinates when approaching the blow-up. 

\begin{figure}[!htb]
\includegraphics[width=0.45\hsize]{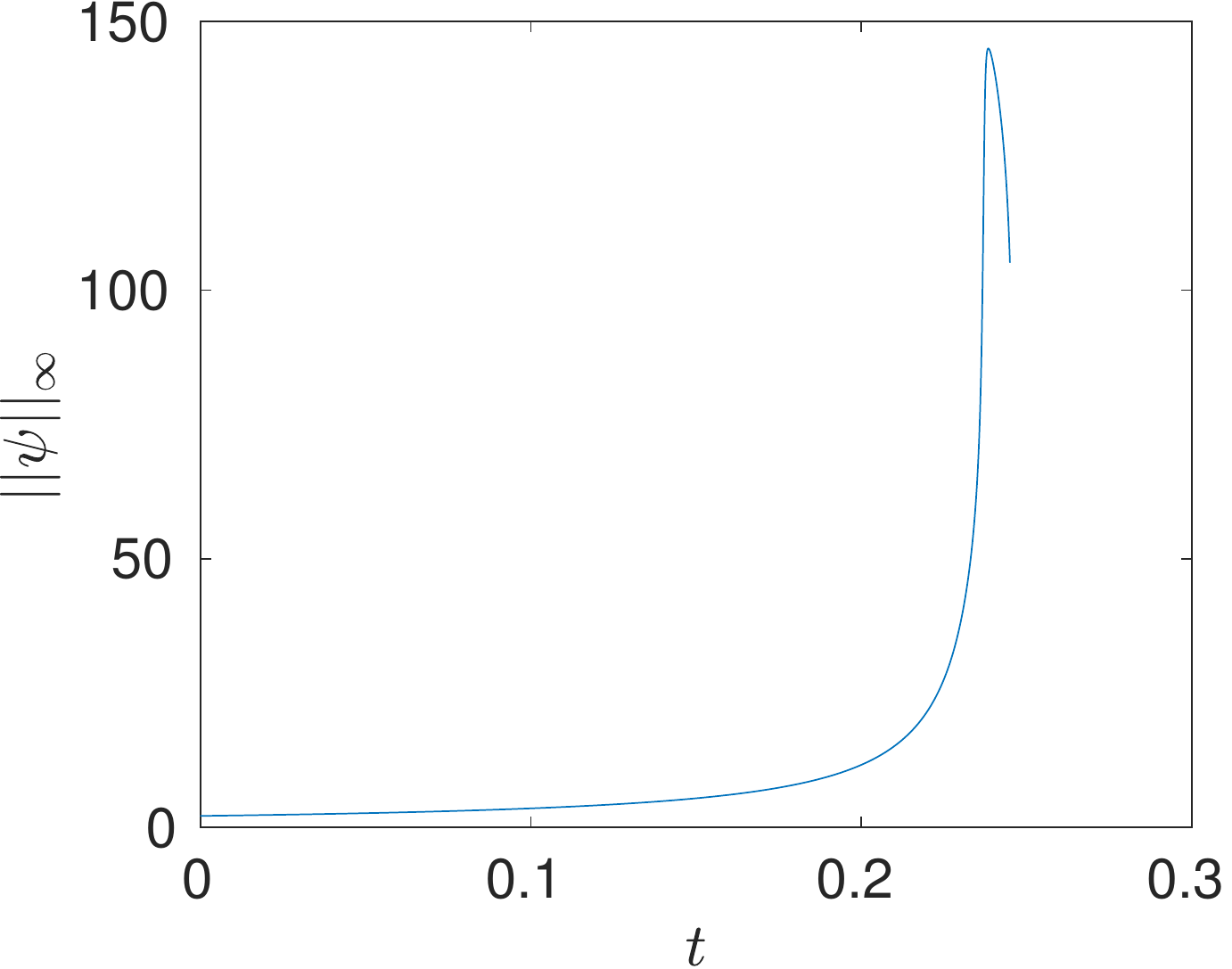}
\includegraphics[width=0.45\hsize]{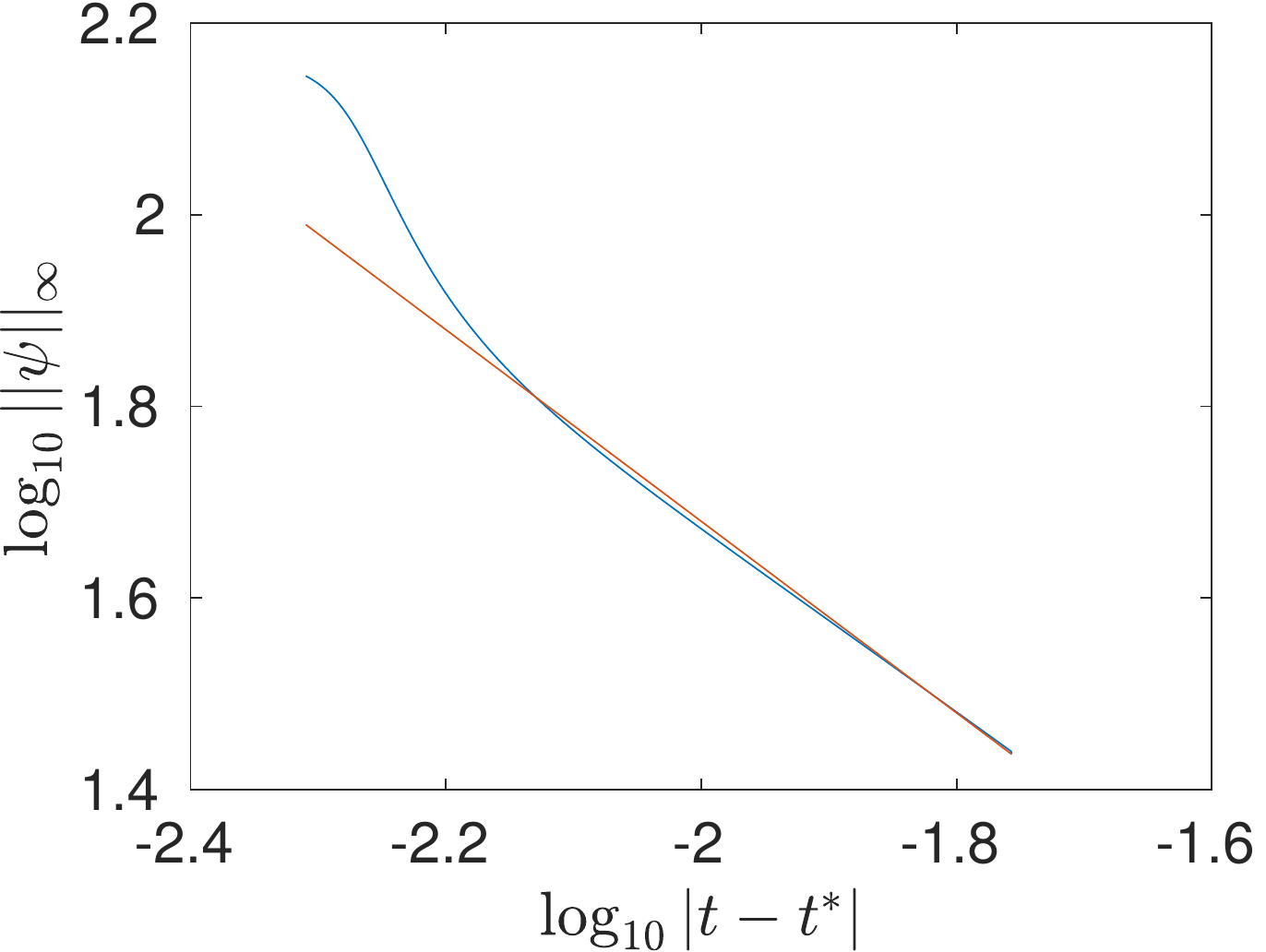}
\caption{Solution to the DS II equation  for Ozawa initial data  
deformed by a small Gaussian $0.1\exp(-x^2 -y^2)$: on  the left the $L^{\infty}$ norm in dependence of time, on the right  a log-log plot of the $L^{\infty}$ 
norm of the solution  near the blow-up together with a straight line 
with slope $-1$.}
\label{figOzdef1b}
\end{figure}

The 
$L^{\infty}$ norm of the solution shown in  Fig.~\ref{figOzdef1b} on 
the left indicates a blow-up before the blow-up time of the 
unperturbed Ozawa solution. 
Fitting the 
$L^{\infty}$ norm for the last 1000 recorded time steps to relation 
(\ref{logfit}) assuming a self similar blow-up provides a blow-up 
time $t^* = 0.2425$  and a blow-up rate $\gamma = - 0.9570$ with a 
fitting residual 0.0033, see Fig.~\ref{figOzdef1b} on the right. Thus 
we get as in the case of the unperturbed Ozawa solution a blow-up 
rate proportional to $t^{*}-t$.

\begin{figure}[!htb]
\includegraphics[width=0.45\hsize]{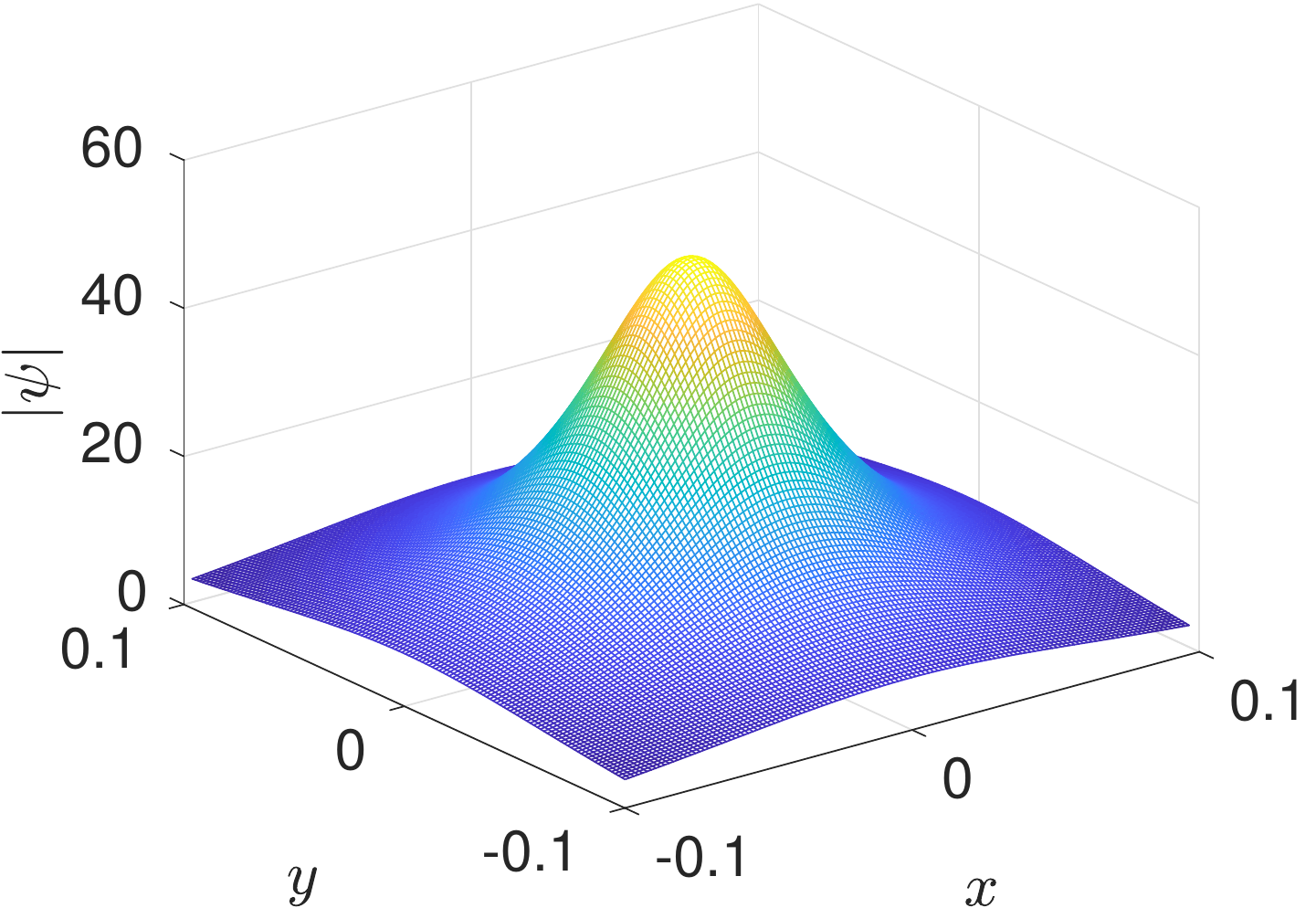}
\includegraphics[width=0.45\hsize]{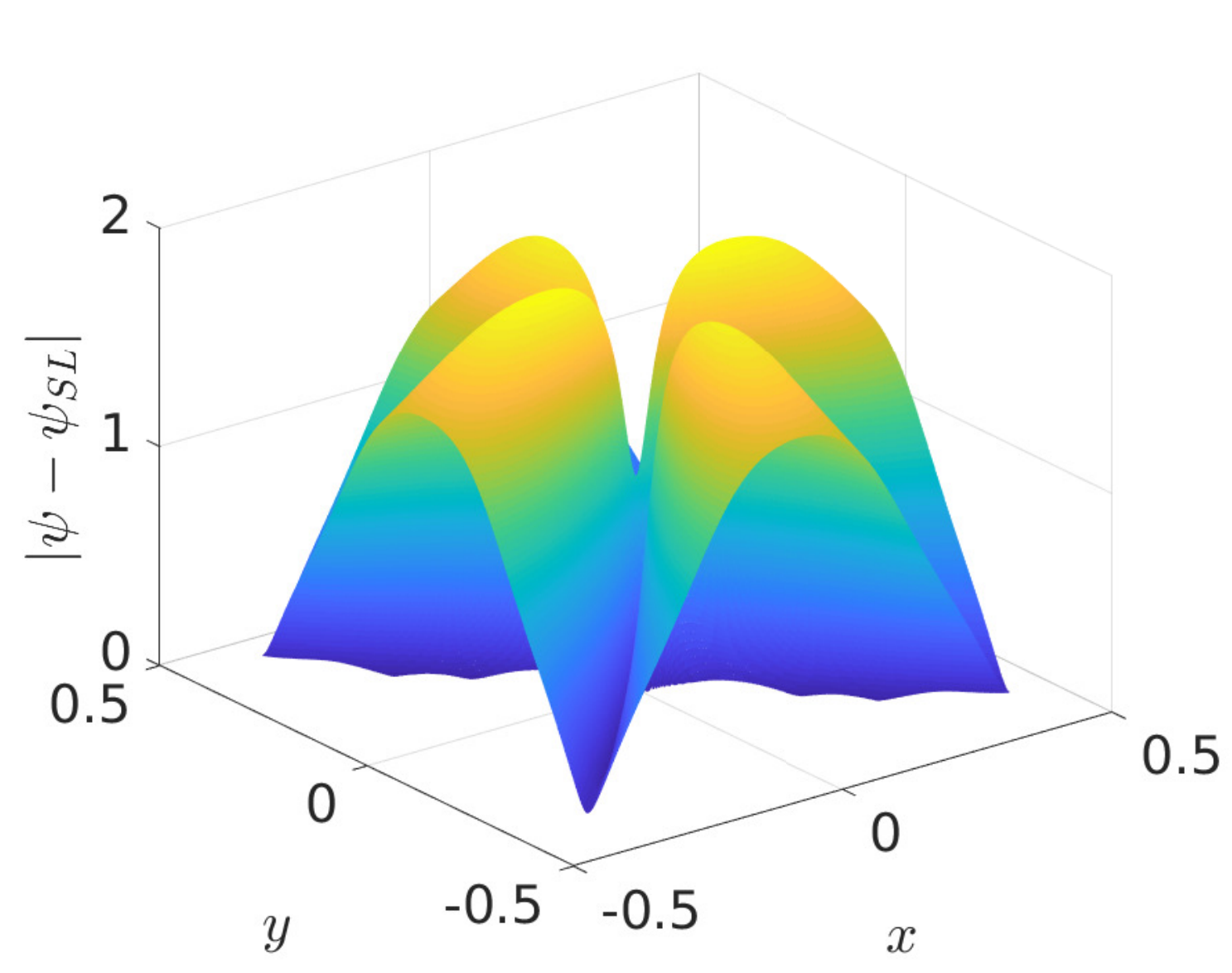}
\caption{
Solution to the DS II equation  for Ozawa initial data  
deformed by a Gaussian $0.1\exp(-x^2 -y^2)$:  blow-up profile at 
$t=0.2326$ on the left and the difference with a rescaled lump 
(according to (\ref{buprofile})) on the right.}
\label{figOzdef1c}
\end{figure}

The blow-up profile at the last recorded time is presented in 
Fig.~\ref{figOzdef1c} on the left. It appears to be locally radially 
symmetric and resembles a lump. Subtracting a lump  rescaled
according to (\ref{buprofile}) ($L$ is just determined via the 
maximum of the solution), we get the figure shown on the right of 
Fig.~\ref{figOzdef1c}. It can be seen that the residual is more than 
an order of magnitude smaller than the maximum of the solution at this 
time. This gives a remarkable agreement with the conjectured 
blow-up profile (\ref{buprofile}) given that it is obtained 
at a time considerably smaller than the actual blow-up time. 

\subsection{Ozawa initial data perturbed by a large Gaussian}
\begin{figure}[!htb]
\includegraphics[width=0.45\hsize]{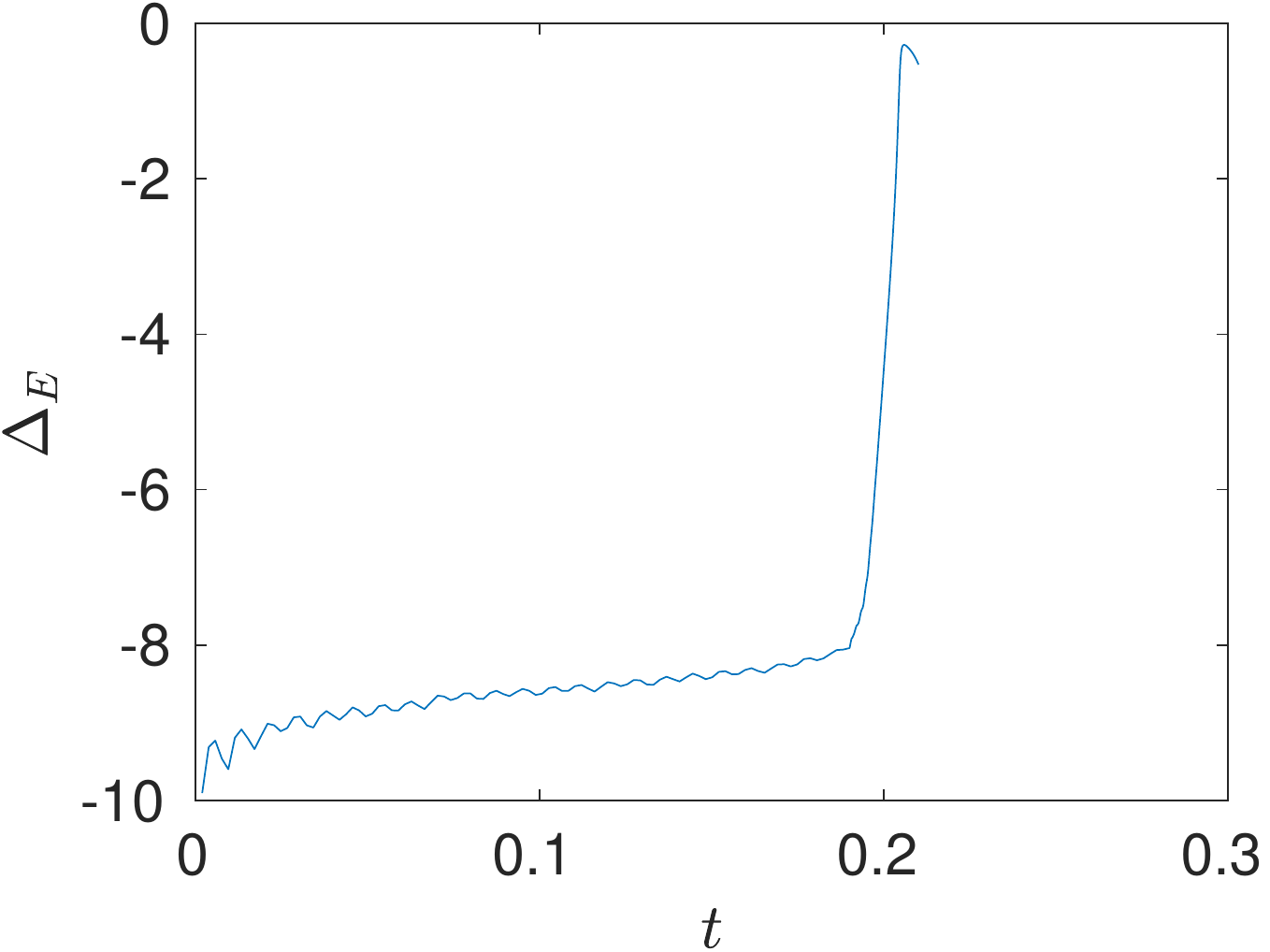}
\includegraphics[width=0.45\hsize]{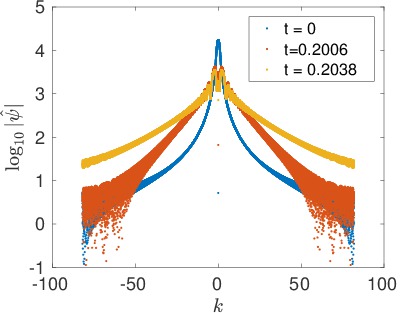}
\caption{Solution to the DS II equation  for Ozawa initial data  
deformed by a Gaussian $0.5\exp(-x^2 -y^2)$: on the left the relative computed energy $\Delta E$, on 
the right the modulus of the Fourier coefficients on the 
$\xi_{1}$-axis  at three 
different times.}
\label{figOzdef5}
\end{figure}

We consider once more Ozawa initial data with parameters $ a = 1$ and 
$b = -4$, but  deform it this time by 
adding a larger Gaussian centered at the origin of the coordinate system  
\begin{equation*}
\psi(x,y,0) = 2\frac{\exp(-i(x^2 - y^2))}{1 + x^2 + y^2} + 0.5\exp{(- x^2 - 
y^2)}.
\end{equation*}
The code is run for 2000 time steps to $t =0.19$, followed by 2000 
time steps $\mathrm{d}t = 10^{-5}$.  The computed relative energy in 
Fig.~\ref{figOzdef5} on the left indicates a loss of resolution in 
time at $t \approx 0.2006$ (step $3060$). In Fig.~\ref{figOzdef5} on 
the right, it can be seen that the Fourier 
coefficients at this time still indicate a spatial resolution to 
better than plotting accuracy. Thus once more, we run out of resolution 
in time first despite a very small time step near the blow-up.

\begin{figure}[!htb]
\includegraphics[width=0.45\hsize]{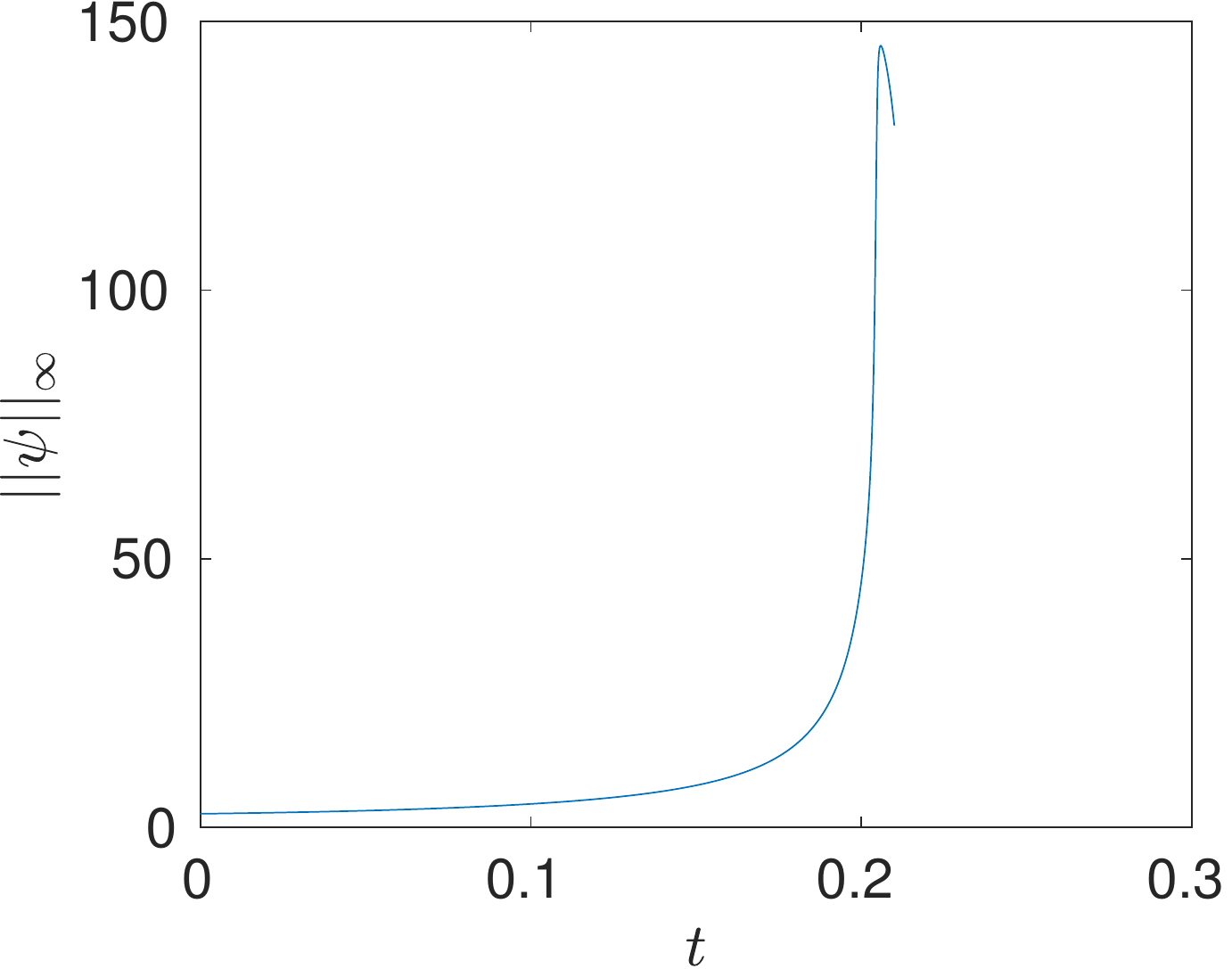}
\includegraphics[width=0.45\hsize]{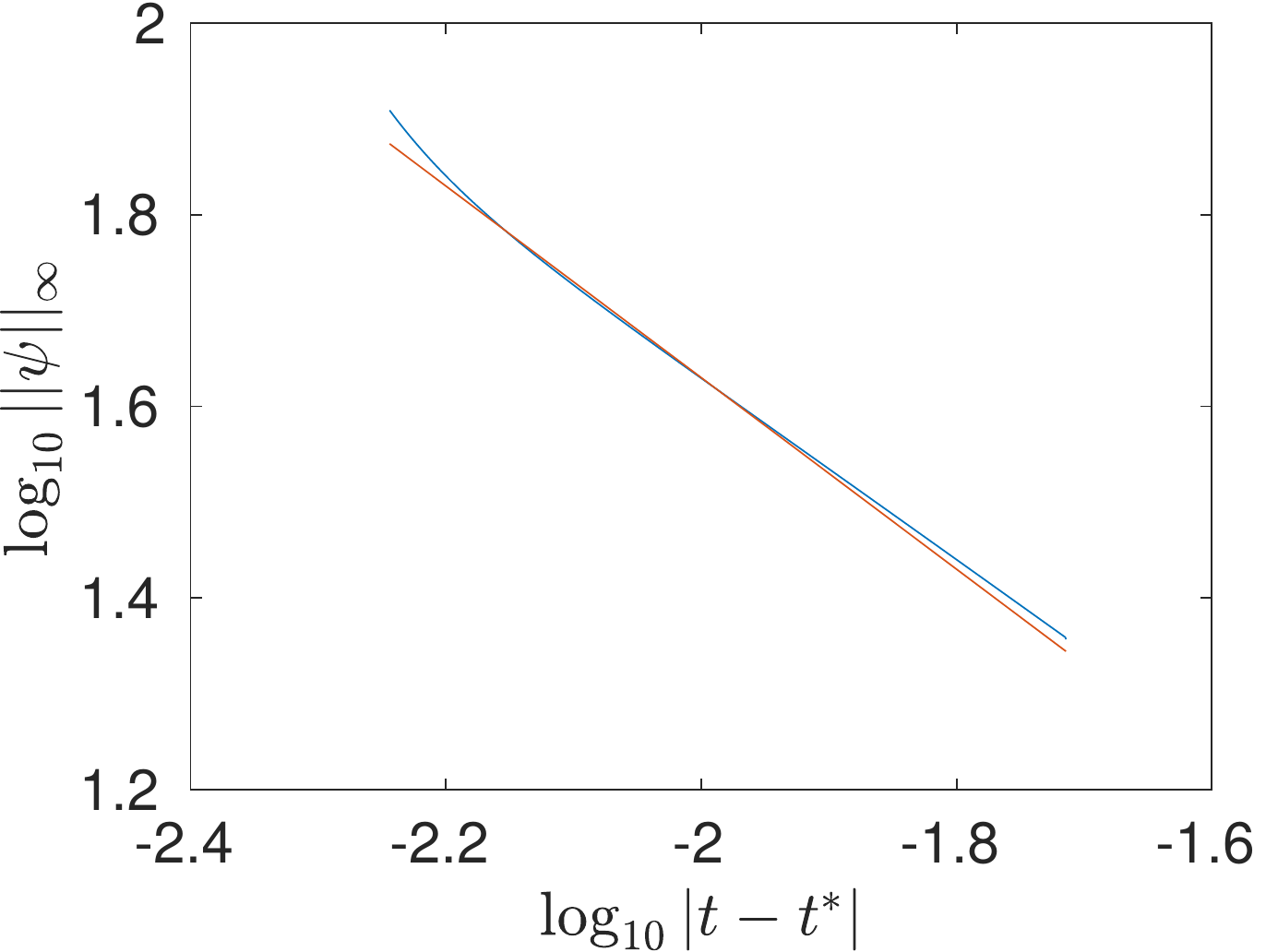}
\caption{Solution to the DS II equation  for Ozawa initial data  
deformed by a Gaussian $0.5\exp(-x^2 -y^2)$: on the left the $L^{\infty}$ norm in dependence of time, on the right  a log-log plot of the $L^{\infty}$ 
norm of the solution  near the blow-up together with a straight line 
with slope $-1$.   }
\label{figOzdef5b}
\end{figure}

As shown in 
Fig.~\ref{figOzdef5b} on the left, the $L^{\infty}$ norm of the 
solution indicates a blow-up at an even earlier time than in 
Fig.~\ref{figOzdef1}. A 
fitting of the $L^{\infty}$ norm 
near blow-up according to \ref{logfit} yields a blow-up time $t^{*} 
=0.2093$ and a blow-up rate $\gamma=-0.9487$ with residual $0.0043$, 
i.e., again the same rate as in the Ozawa solution.

\begin{figure}[!htb]
\includegraphics[width=0.45\hsize]{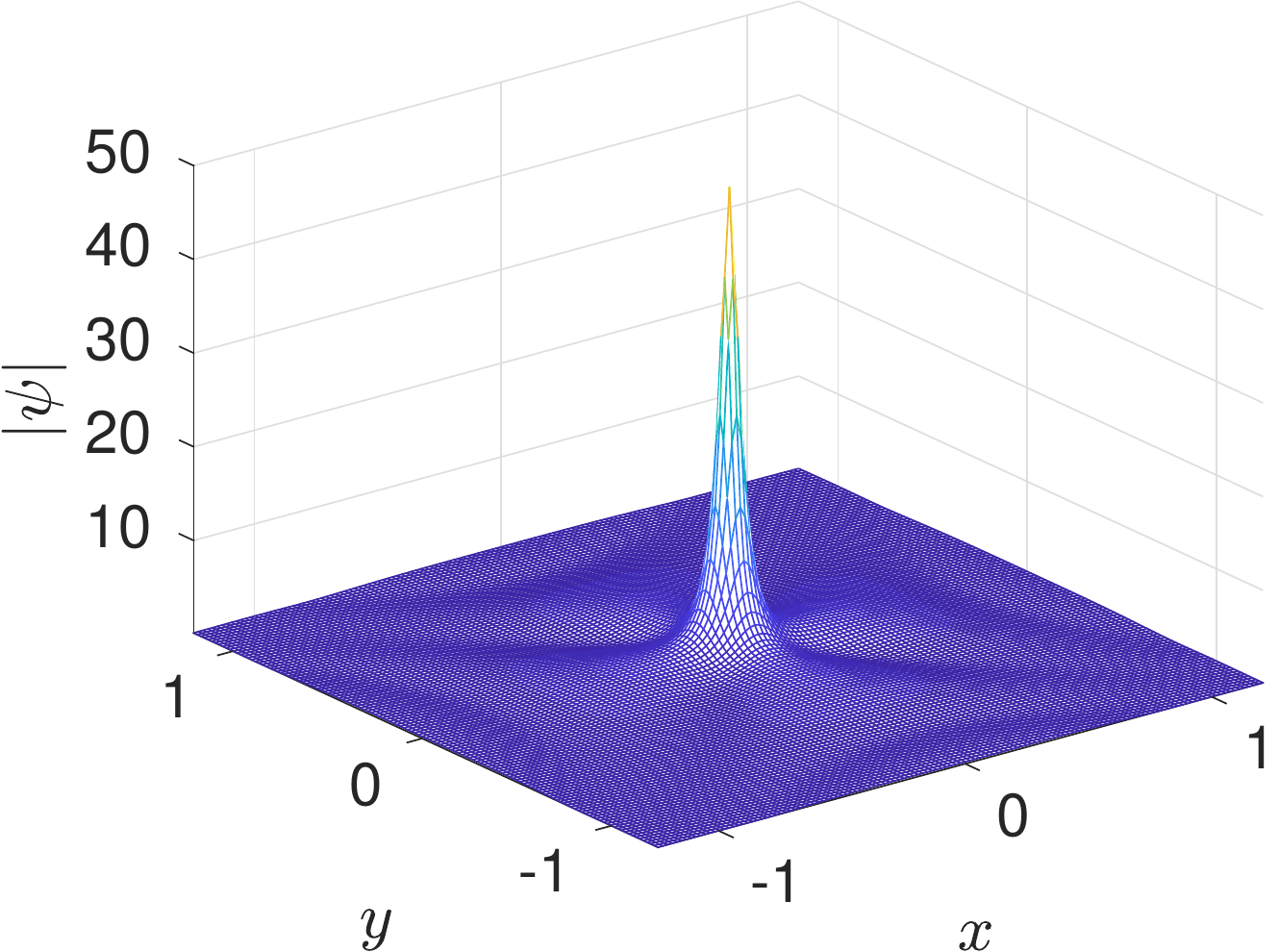}
\includegraphics[width=0.45\hsize]{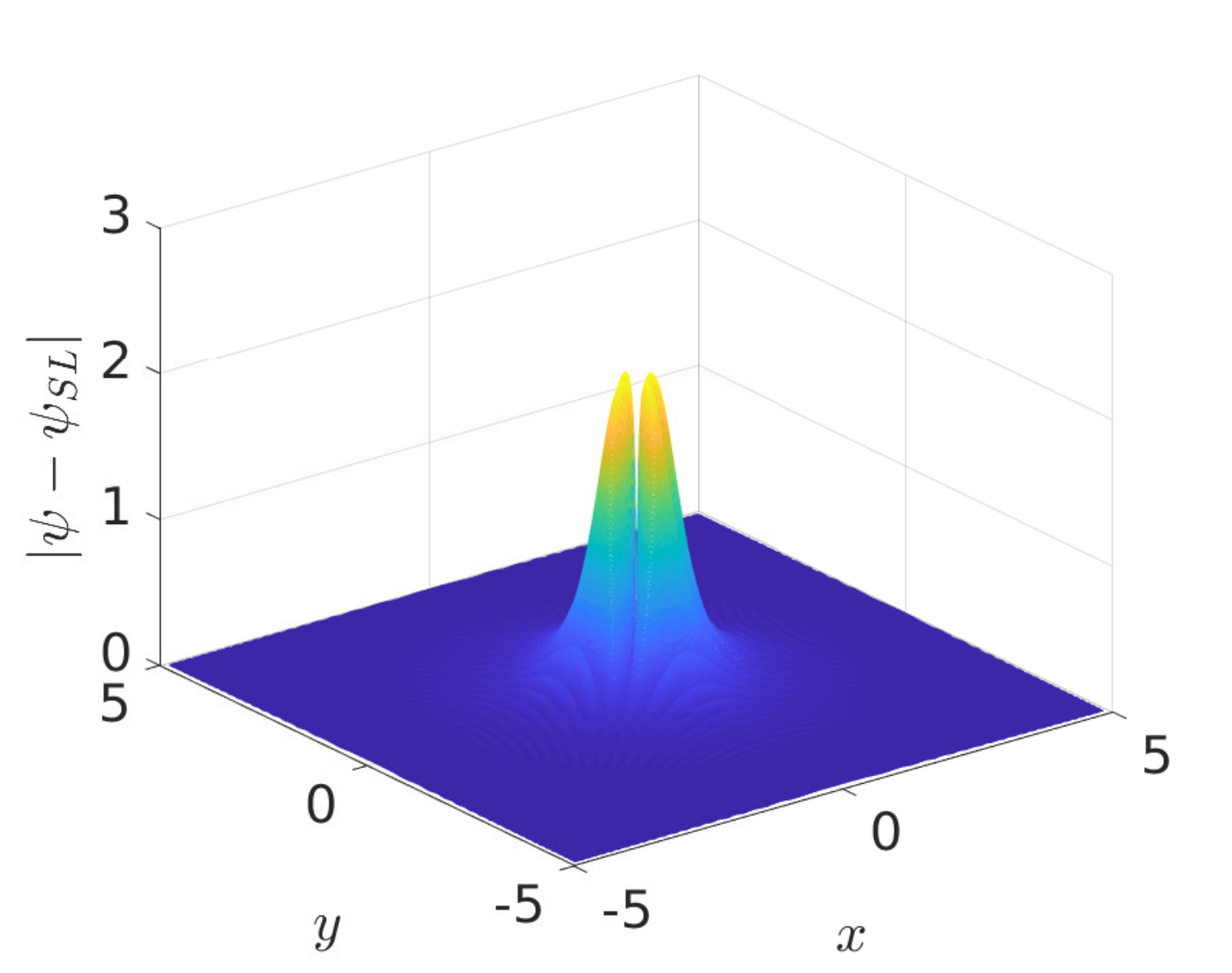}
\caption{Solution to the DS II equation  for Ozawa initial data  
deformed by a Gaussian $0.5\exp(-x^2 -y^2)$:  blow-up profile at 
$t=0.2006$ on the left and the difference with a rescaled lump 
(according to (\ref{buprofile})) on the right.}
\label{figOzdef5c}
\end{figure}

The blow-up profile at the last recorded time is presented in 
Fig.~\ref{figOzdef5c} on the left. Once more it seems to be  radially 
symmetric in the vicinity of the maximum. We show the solution here 
on larger scales than in Fig.~\ref{figOzdef1c} to allow an overview 
of the solution, a close-up would look as in Fig.~\ref{figOzdef1c}.   
If we subtract as there a lump  rescaled
according to (\ref{buprofile}), we get the figure shown on the right of 
Fig.~\ref{figOzdef5c}. Once more the residual is more than 
an order of magnitude smaller than the maximum of the solution at this 
time, which shows an excellent agreement with the asymptotic blow-up 
profile. 

\subsection{Perturbed lump solution}
\begin{figure}[!htb]
\includegraphics[width=0.45\hsize]{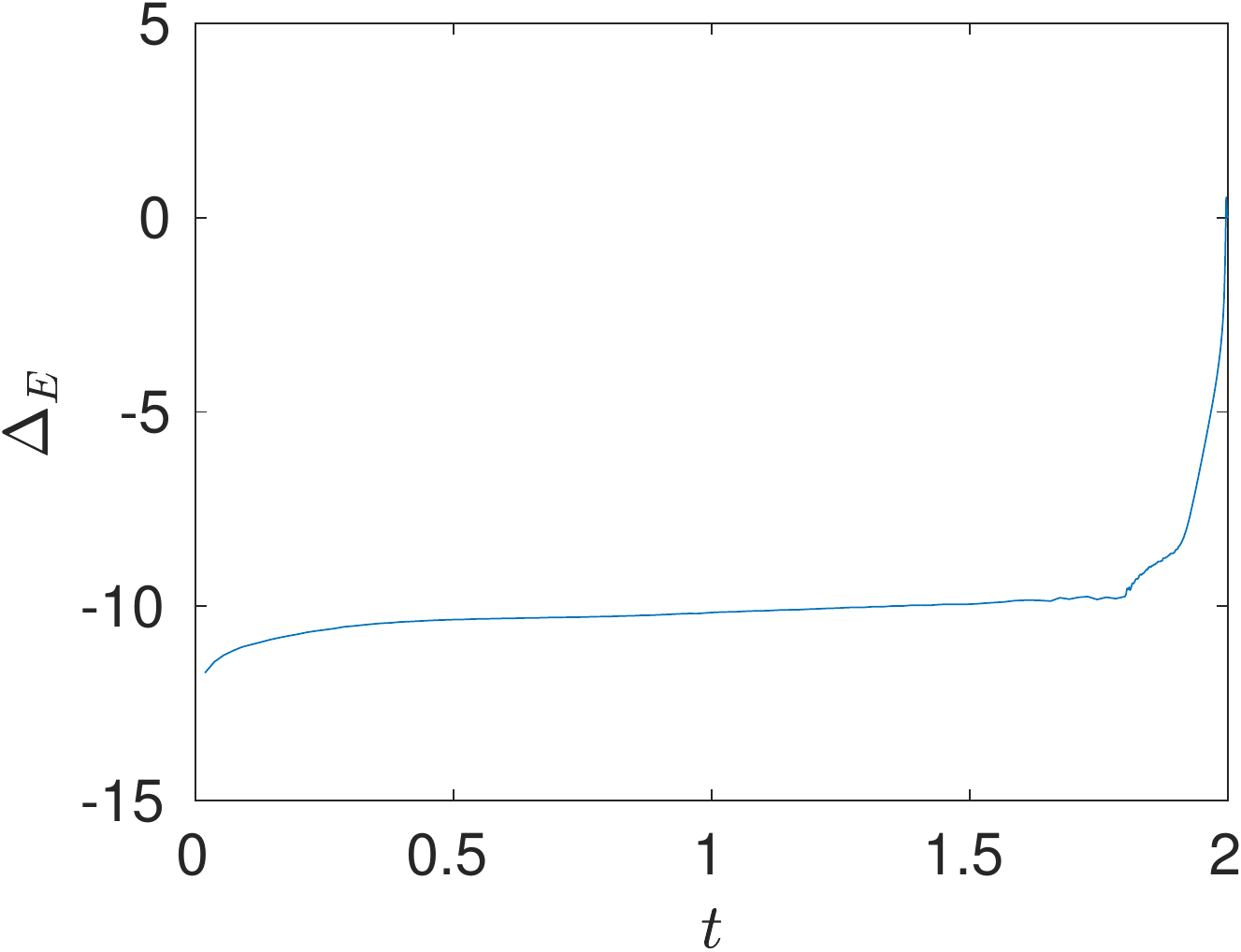}
\includegraphics[width=0.45\hsize]{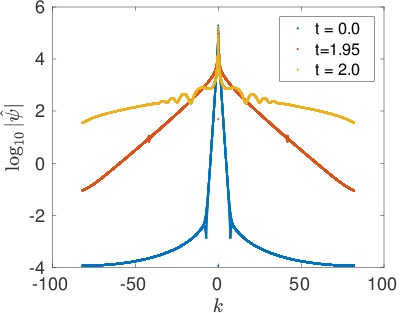}
\caption{Solution to the focusing DS II equation for deformed lump 
initial data  $ \psi(x,y,0) = 2.2/(1 + x^2 +y^2)$: on the left the relative computed energy $\Delta E$, on 
the right the modulus of the Fourier coefficients on the $\xi_{1}$-axis at three 
different times.}
\label{figlump}
\end{figure}

We consider lump initial data and multiply it by a factor $1.1$, i.e., 
we take the initial condition 
\begin{equation*}
\psi(x,y,0) = \frac{2.2}{1 + x^2 + y^2}.
\end{equation*}
The code is run for 2000 steps to $t =1.9$, followed by 2000 time steps 
 $\mathrm{d}t = 10^{-4}$. The time step $\mathrm{d} t$ is chosen so 
that $\mathrm{d} t/t^*$ is of the same order of magnitude in all experiments. A jump in the energy for 
$t\sim2$ as shown in Fig.~\ref{figlump} on the left indicates again a loss of 
resolution in $t$, whereas the Fourier coefficients on the right of 
Fig.~\ref{figlump} also degrade for $t\sim2$, again slightly later than 
the time where the energy has a jump. 

\begin{figure}[!htb]
\includegraphics[width=0.45\hsize]{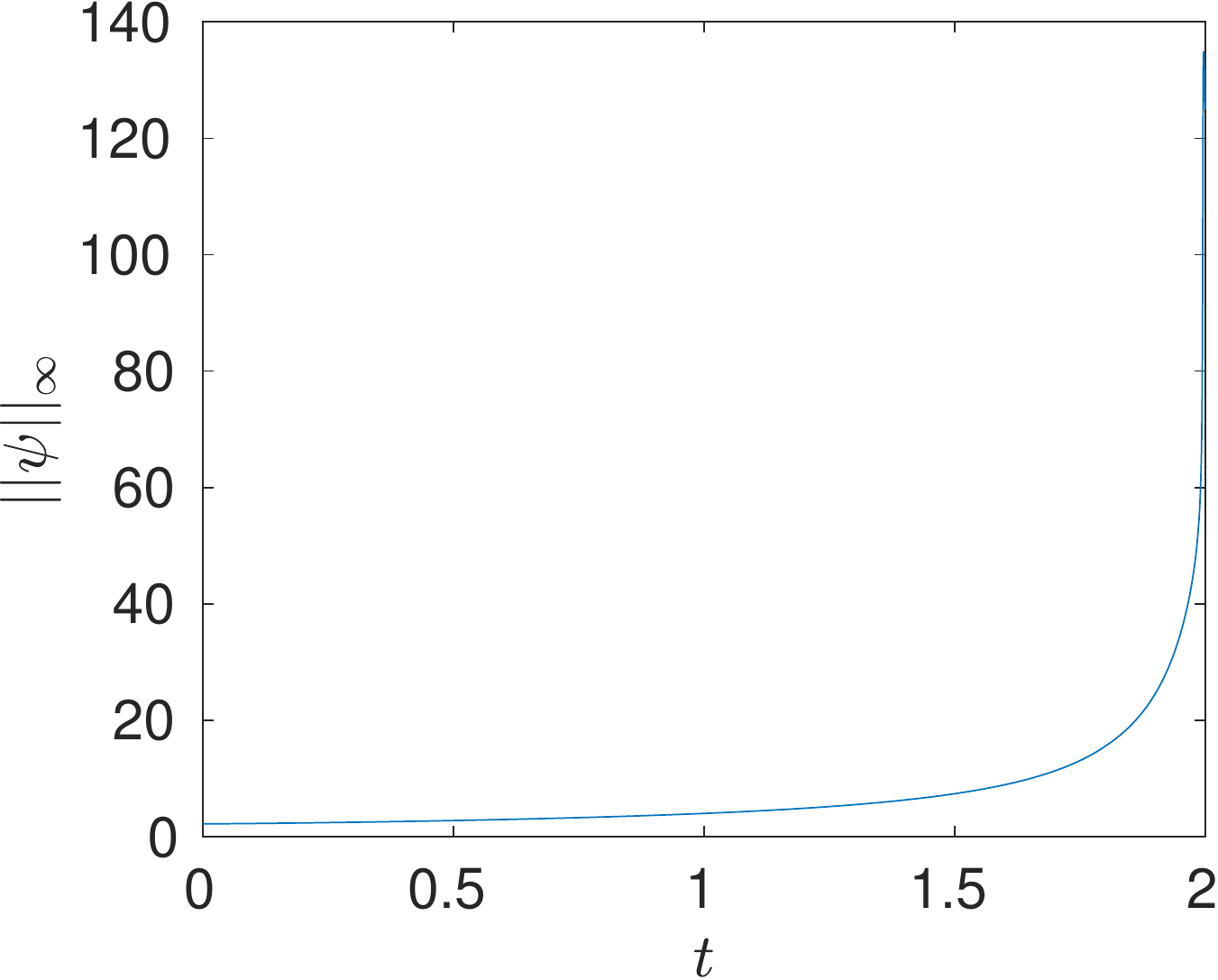}
\includegraphics[width=0.45\hsize]{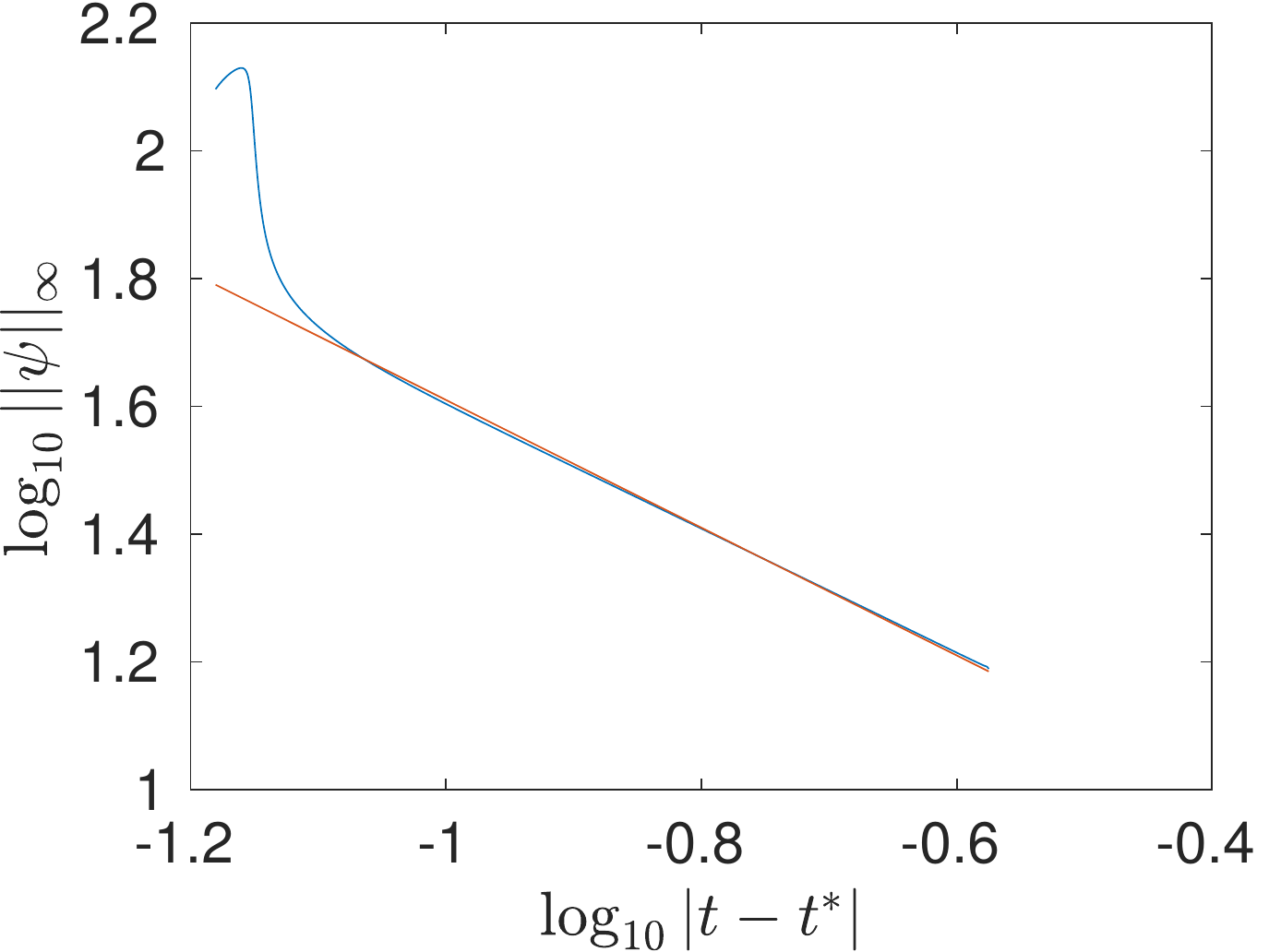}\\
\caption{Solution to the focusing DS II equation for deformed lump 
initial data  $ \psi(x,y,0) = 2.2/(1 + x^2 +y^2)$: on the left the $L^{\infty}$ norm in dependence of time, on the right  a log-log plot of the $L^{\infty}$ 
norm of the solution  near the blow-up together with a straight line 
with slope $-1$. }
\label{figlumpb}
\end{figure}

The $L^{\infty}$ norm of the solution in Fig.~\ref{figlumpb} on the 
left is once more indicating a blow-up for $t>2$. The fit of the 
$L^{\infty}$ norm for the last 1000 time steps to (\ref{logfit}) 
yields a blow-up time $t^* = 2.066$ with blow-up rate $\gamma = 
-0.9733$ with a residual of 0.0043 as can be seen on the right of 
Fig.~\ref{figlumpb}.

\begin{figure}[!htb]
\includegraphics[width=0.45\hsize]{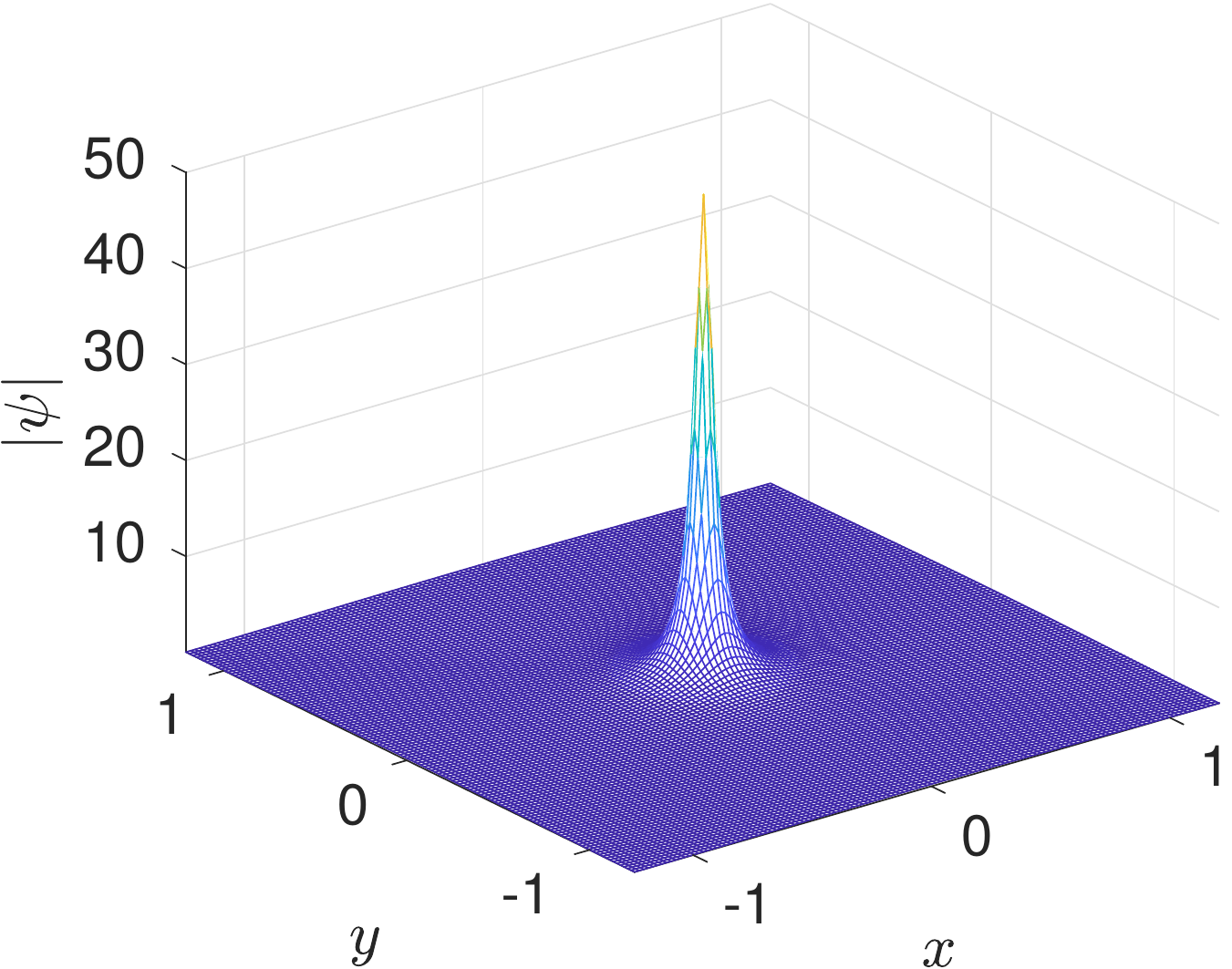}
\includegraphics[width=0.45\hsize]{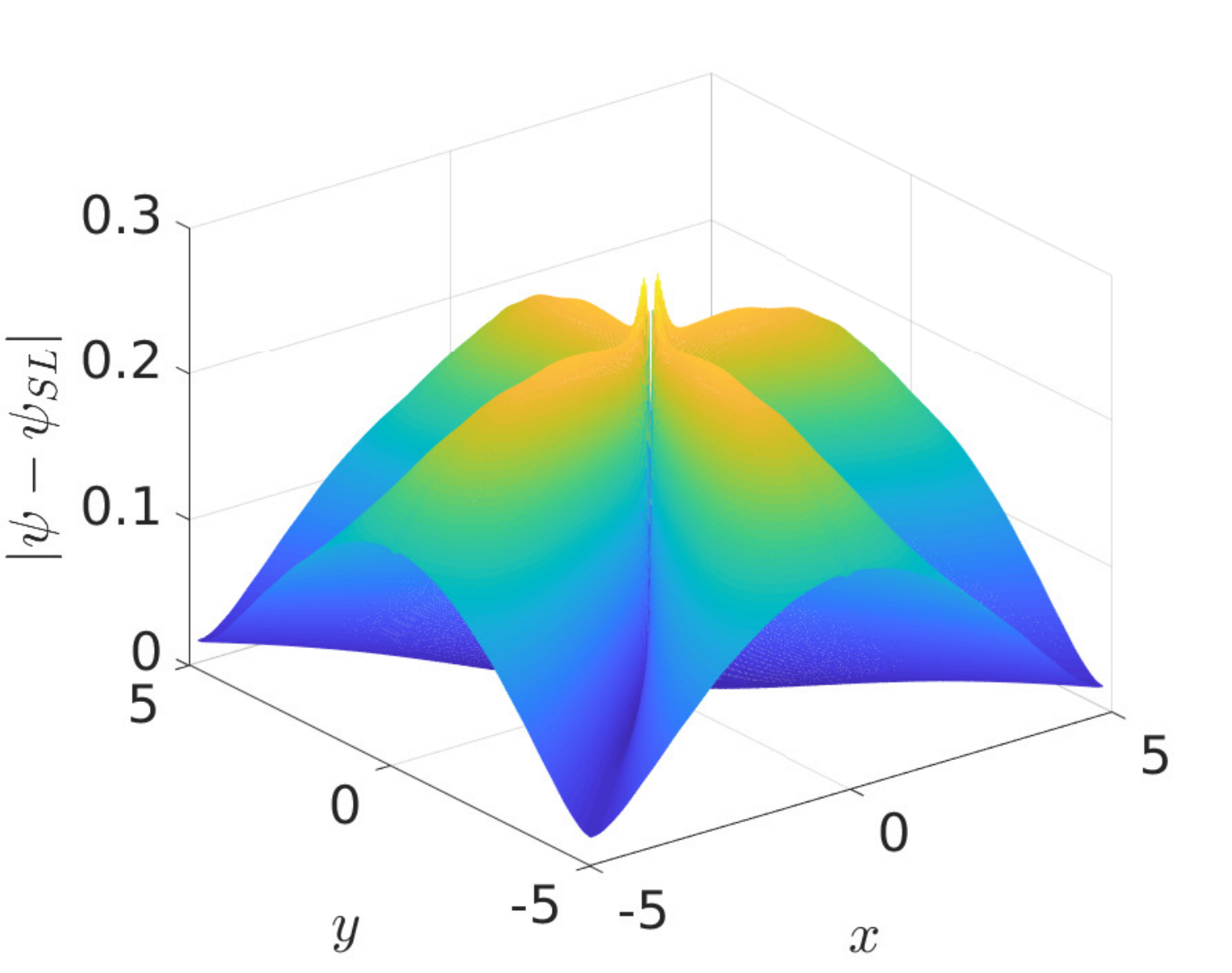}
\caption{Solution to the focusing DS II equation for deformed lump 
initial data  $ \psi(x,y,0) = 2.2/(1 + x^2 +y^2)$: blow-up profile at 
$t=1.983$ on the left and the difference with a rescaled lump 
(according to (\ref{buprofile})) on the right.}
\label{figlumpc}
\end{figure}

The blow-up profile for $t=2$ can be seen in Fig.~\ref{figlumpc} on 
the left. Here the agreement with a dynamically rescaled lump 
(\ref{buprofile}) is even 
more striking than in the Ozawa examples above. As can be seen on the 
right of Fig.~\ref{figlumpc}, the residual of the numerical solution 
and the rescaled lump is almost two orders of magnitude smaller than 
the maximum of the solution. 

\subsection{Gaussian initial data}
If one is interested in the solution to the DS II equation for 
localized initial data varying on length scales of order 
$1/\epsilon$, and this for times of order $1/\epsilon$ with 
$\epsilon\ll1$, a way to treat this is a change of coordinates 
$x\mapsto \epsilon x$, $y\mapsto \epsilon y$, and $t\mapsto \epsilon 
t$. This leads for (\ref{DSII}) to the equation (in an abuse of notation 
we use the same symbols as before)
\begin{equation}
    \label{DSIIe}
\begin{array}{ccc}
i\epsilon
\partial_{t}\psi+\epsilon^{2}\partial_{xx}\psi-\epsilon^{2}\partial_{yy}\psi-2\left(\Phi+\left|\psi\right|^{2}\right)\psi & = & 0,
\\
\partial_{xx}\Phi+\partial_{yy}\Phi+2\partial_{xx}\left|\psi\right|^{2} & = & 0.
\end{array}
\end{equation}
Thus it is possible to consider a family of equations (depending on 
the parameter $\epsilon$) for $\epsilon$-dependent initial data, 
which allows to study the semiclassical limit of DS II, see for 
instance
\cite{KR2,AKMM}.

For equation (\ref{DSIIe}) with $\epsilon=0.1$, we consider Gaussian initial data
\begin{equation*}
\psi(x,y,0) = \exp(- x^2 - y^2)
\end{equation*}
on a domain with $D=2$. The computed energy jumps at $t = 0.2954$ 
(step 3540), whereas the fitting of the Fourier coefficients 
according to (\ref{fourasymp}) indicates a singularity at a distance 
from the real axis smaller than the spatial resolution at $t=0.2956$ 
(step 3557). Thus the loss of resolution  in time and in space 
happens at roughly the same time in this example.  
\begin{figure}[!htb]
\includegraphics[width=0.45\hsize]{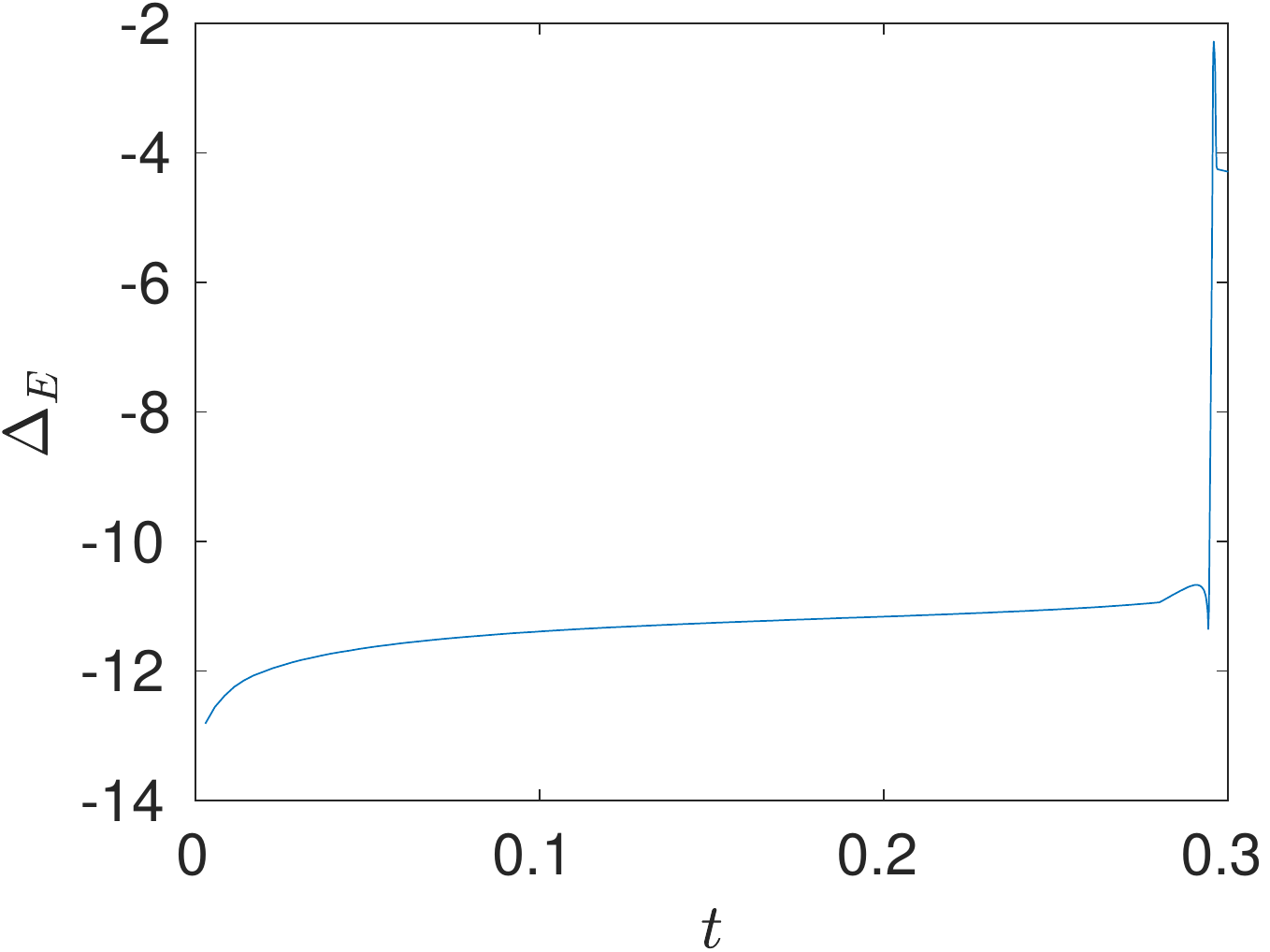}
\includegraphics[width=0.45\hsize]{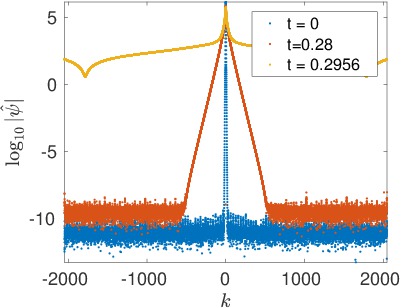}
\caption{Solution to the focusing DS II equation (\ref{DSIIe}) for  
Gaussian initial data  initial data  $\psi(x,y,0) =\exp(-x^2 -y^2)$ and 
$\epsilon=0.1$: on the left the relative computed energy $\Delta E$, on 
the right the modulus of the Fourier coefficients on the $\xi_{1}$-axis at three 
different times.}
\label{figGauss}
\end{figure}

\begin{figure}[!htb]
\includegraphics[width=0.45\hsize]{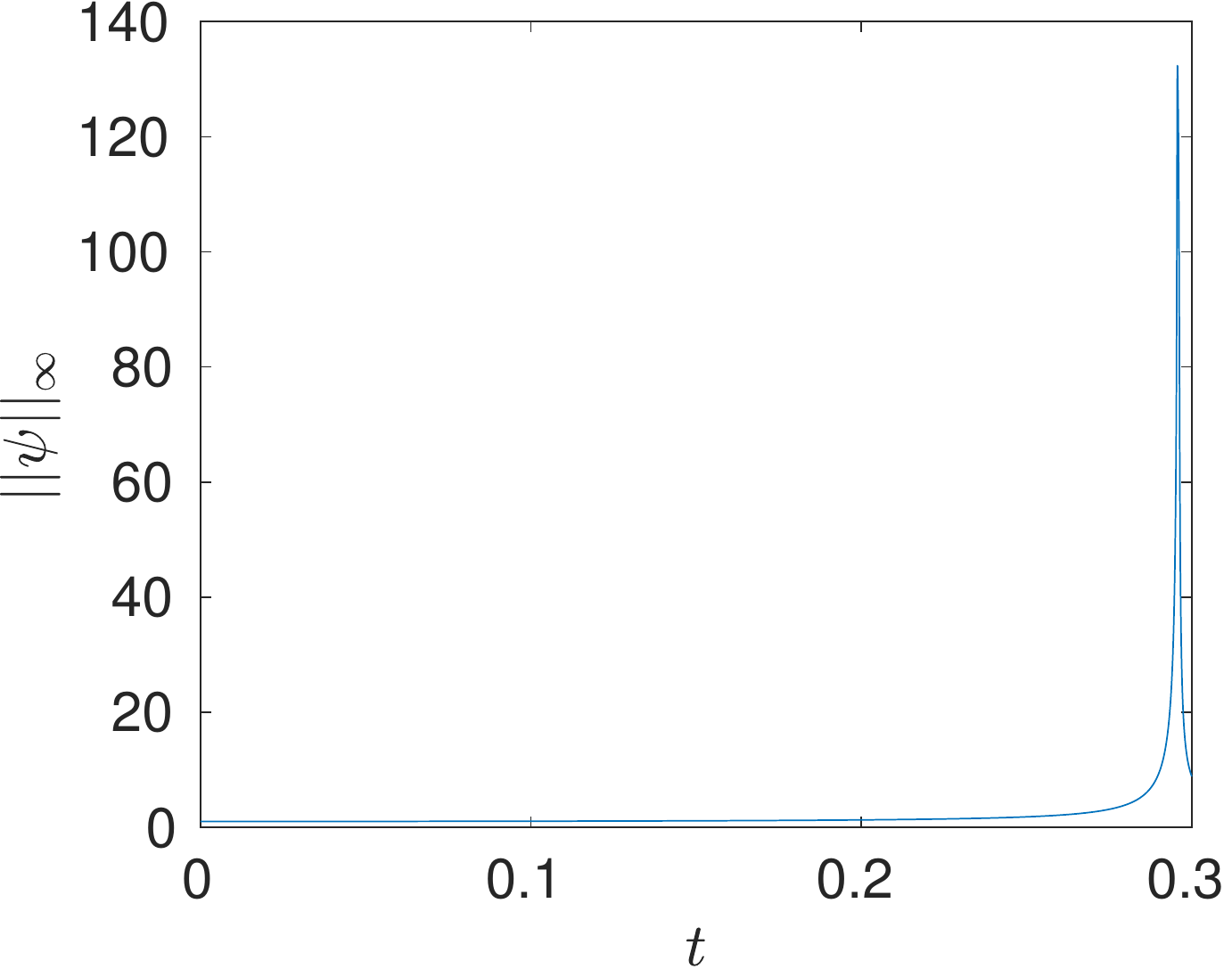}
\includegraphics[width=0.45\hsize]{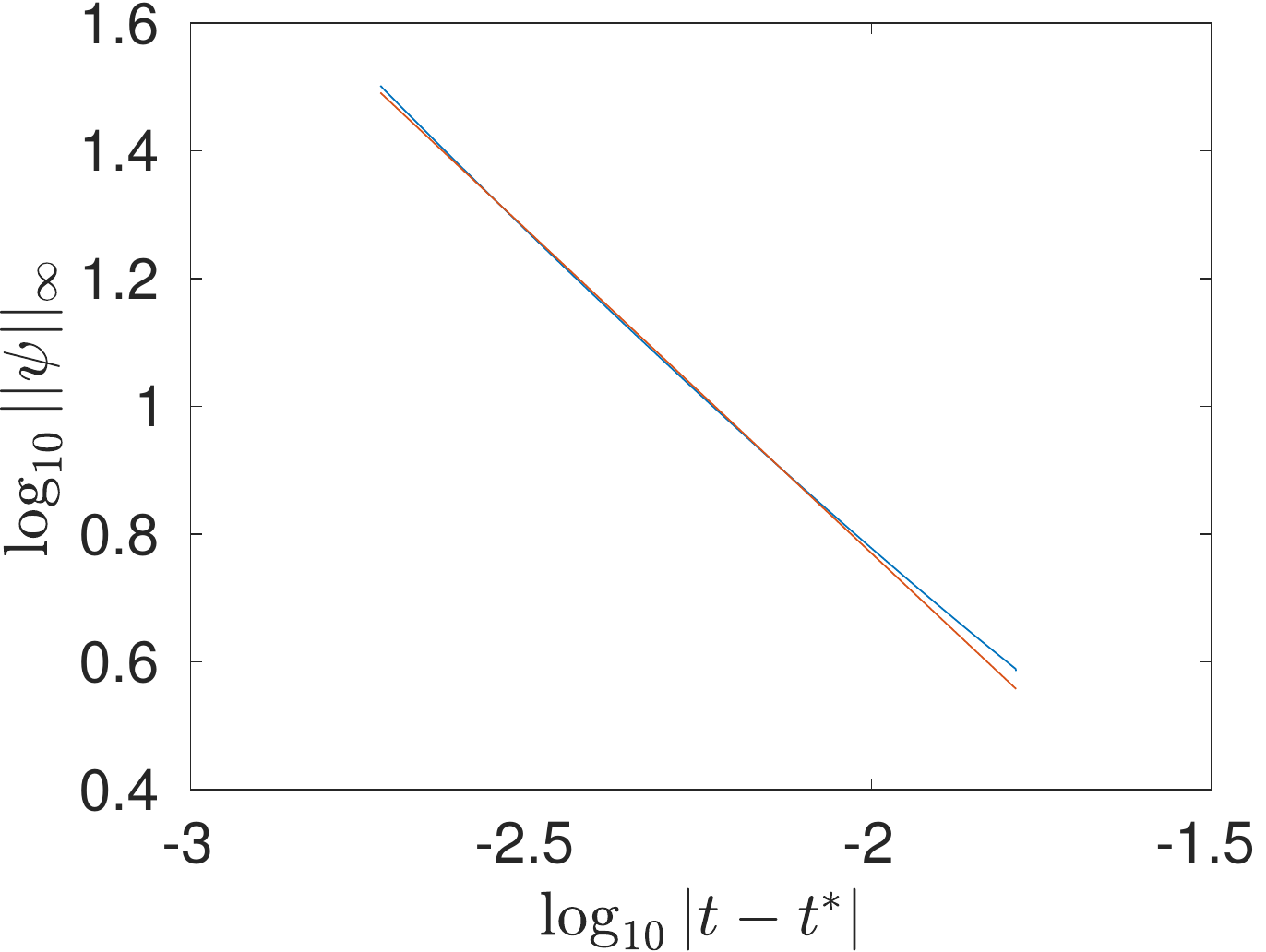}
\caption{Solution to the focusing DS II equation (\ref{DSIIe}) for  
Gaussian initial data  initial data  $\psi(x,y,0) =\exp(-x^2 -y^2)$ and 
$\epsilon=0.1$: on the left the $L^{\infty}$ norm in dependence of time, on the right  a log-log plot of the $L^{\infty}$ 
norm of the solution  near the blow-up together with a straight line 
with slope $-1$.}
\label{figGaussb}
\end{figure}

The $L^{\infty}$ norm on the left of Fig.~\ref{figGaussb} indicates a 
rapid blow-up. Fitting the norm according to (\ref{logfit}) gives a 
blow-up time $t^{*} = 0.2965$ and a blow-up rate $\gamma=-0.9278$. In 
this case, due to the exponential decay of the solution in the 
spatial directions we are able to use a much smaller computational 
domain and thus to work with a much higher spatial resolution as 
before. This allowed to get much closer to the actual blow-up time 
than in the algebraically decaying examples above. 

\begin{figure}[!htb]
\includegraphics[width=0.45\hsize]{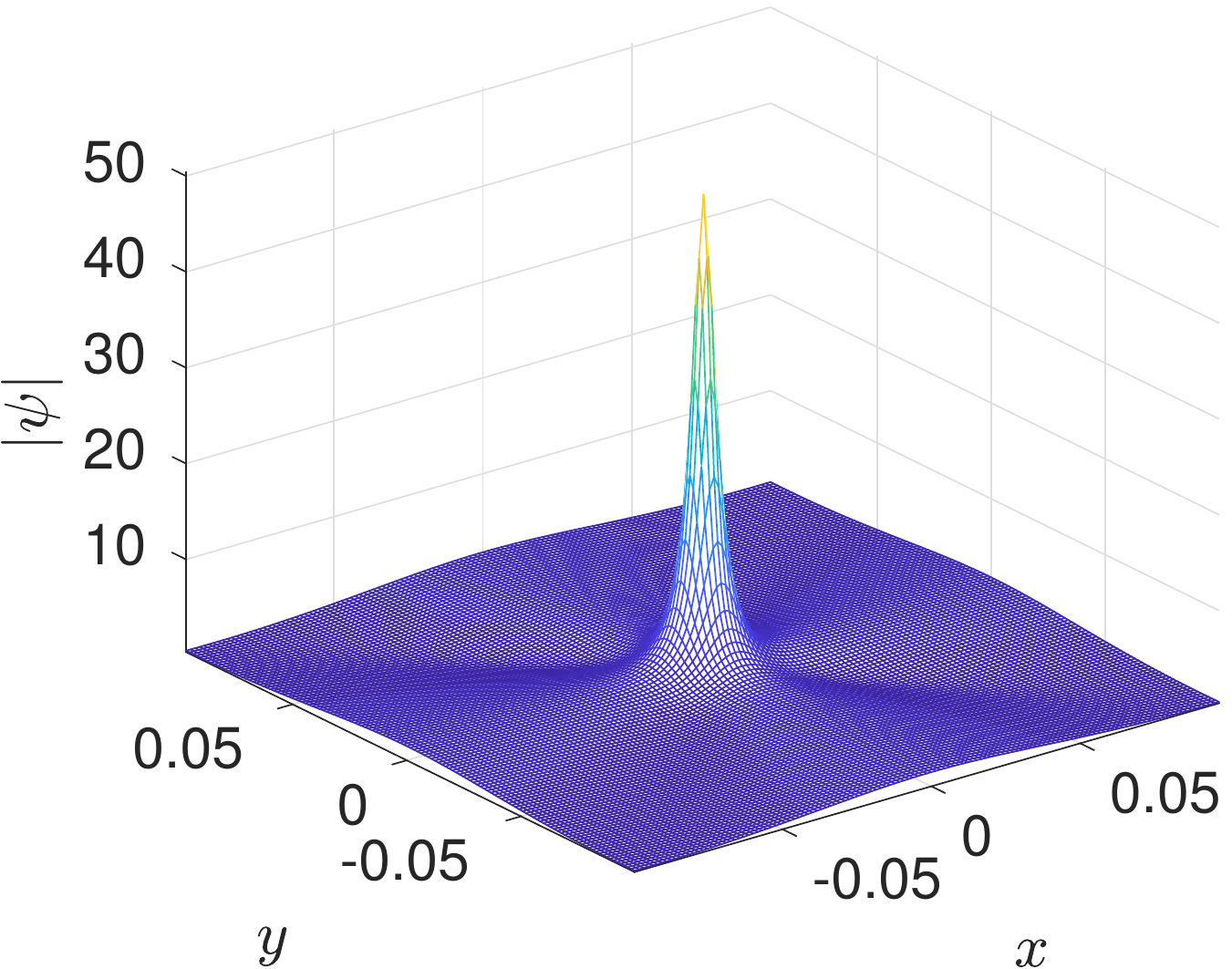}
\includegraphics[width=0.45\hsize]{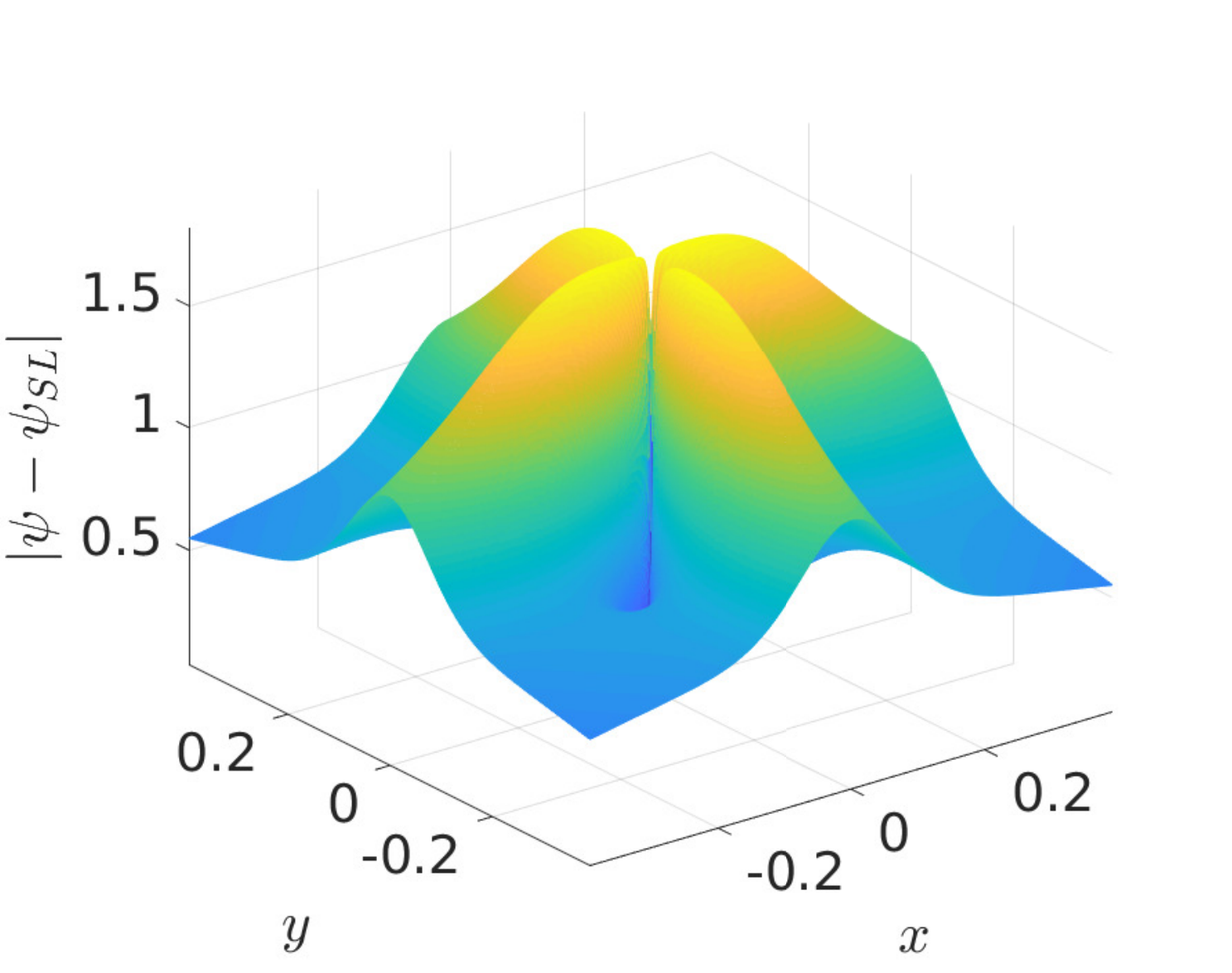}
\caption{Solution to the focusing DS II equation (\ref{DSIIe}) for  
Gaussian initial data  initial data  $\psi(x,y,0) =\exp(-x^2 -y^2)$ and 
$\epsilon=0.1$:  blow-up profile at 
$t=0.2950$ on the left and the difference with a rescaled lump 
(according to (\ref{buprofile})) on the right.}
\label{figGaussc}
\end{figure}

The blow-up profile can be seen in Fig.~\ref{figGaussc} on the left. 
Subtracting a rescaled lump (according to (\ref{buprofile}), we 
obtain the residual on the right of the same figure. This indicates 
that also in this case, the blow-up profile appears to be given by a 
dynamically rescaled lump. 

\section{Outlook}
In this paper we have studied blow-up phenomena in the focusing DS II 
equation for asymptotically localized initial data with a single 
maximum of the modulus. It was shown that initial data both with algebraic and with 
exponential decay can lead to a blow-up if the $L^{\infty}$ norm of 
the initial data is sufficiently large. Already in \cite{KR2013a} and 
\cite{KS} it was shown that blow-up occurs only if the initial data 
are locally radially symmetric near the maximum. It would be 
interesting to study this symmetry question in more detail since it 
appears to be related to the lump appearing as a blow-up profile in a 
self-similar blow-up.  This question would be in particular 
intriguing in the context of several humps which could blow up at 
the same time if the initial data are carefully chosen. 

To study such questions, higher resolution than available for the 
present study will be needed. A possible way to address this would be 
to take a multi-domain spectral approach as in \cite{BK} and 
references therein which would 
be especially beneficial in the case of an algebraic decay of the 
solution towards infinity. It would also allow to allocate more 
resolution near a blow-up by introducing adapted domains. In contrast 
to \cite{BK} where a Laplace operator in the NLS equation 
was studied, DS II equations have a d'Alembert operator for which a 
different compactification at infinity will have to be introduced. 
This is will be the topic of future work. 

\begin{merci} This work was partially supported by the program PARI and the FEDER 
2016 and 2017. We thank J.-C. Saut for helpful discussions and hints. 
\end{merci}

\end{document}